\documentclass[11pt,letterpaper]{article}
\usepackage{amsfonts}
\usepackage{amsmath,amssymb,color}
\usepackage{graphicx}
\usepackage[singlespacing]{setspace}

\def\rr{\mbox{$\mathbf r$}}

\def\xx{\mbox{$\mathbf x$}}

\newcommand{\yy}{\mbox{$\mathbf y$}}

\newtheorem{theorem}{Theorem}

\newtheorem{definition}{Definition}
\newtheorem{example}{Example}

\newtheorem{lemma}{Lemma}

\newtheorem{proposition}{Proposition}
\newtheorem{remark}{Remark}

\newenvironment{proof}[1][Proof]{\noindent\textbf{#1.} }{\ \rule{0.5em}{0.5em}}

\begin{document}

\title{On the Stationary Distribution of Iterative Imputations}

\author{Jingchen Liu, Andrew Gelman, Jennifer Hill, and Yu-Sung Su\\ \\
Columbia University, Columbia University, New York University, \\  Tsinghua University}
\date{\today}
\maketitle
\begin{abstract}
Iterative imputation, in which variables are imputed one at a time each given a model predicting from all the others, is a popular technique that can be convenient and flexible, as it replaces a potentially difficult multivariate modeling problem with relatively simple univariate regressions. In this paper, we begin to characterize the stationary distributions of iterative imputations and their statistical properties. More precisely, when the conditional models are compatible (defined in the text), we give a set of sufficient conditions under which the imputation distribution converges in total variation to the posterior distribution of a Bayesian model. When the conditional models are incompatible but are valid, we show that the combined imputation estimator is consistent.
\end{abstract}


\baselineskip=18pt
\section{Introduction}
Iterative imputation is a widely used approach for imputing multivariate missing data.  The procedure starts by randomly imputing missing values using some simple stochastic algorithm. Missing values are then imputed one variable at a time, each conditionally on all the others using a model fit to the current iteration of the completed data.  The variables are looped through until approximate convergence (as measured, for example, by the mixing of multiple chains). 

Iterative imputation can be an appealing way to express uncertainty about missing data.  There is no need to explicitly construct a joint multivariate model of all types of variables: continuous, ordinal, categorical, and so forth. Instead, one need only specify a sequence of families of conditional models such as linear regression, logistic regression, and other standard and already programmed forms. The distribution of the resulting imputations is implicitly defined as the invariant (stationary) distribution of the Markov chain corresponding to the iterative fitting and imputation process.

Iterative, or chained, imputation is convenient and flexible and has been implemented in various ways in several statistical software packages, including {\tt mice} \cite{MICE10} and {\tt mi} \cite{MI10} in \texttt{R}, {\tt IVEware} \cite{IVEware10} in \texttt{SAS}, and {\tt ice} in \texttt{Stata} \cite{ICE04,ICE05}. The popularity of these programs suggests that the resulting imputations are believed to be of practical value.  However, the theoretical properties of iterative imputation algorithms are not well understood.  Even if, as we would prefer, the fitting of each imputation model and the imputations themselves are performed using conditional Bayesian inference, the stationary distribution of the algorithm (if it exists) does not in general correspond to Bayesian inference on any specified multivariate distribution.

Key questions are:  (1) Under what conditions does the algorithm converge to a stationary distribution?  (2) What statistical properties does the procedure admit given that a unique stationary distribution exists?

Regarding the first question, researchers have long known that
the Markov chain may be non-recurrent (``blowing up'' to infinity or drifting like a nonstationary random walk), even if each of the conditional models is fitted using a proper prior distribution. 

In this paper, we focus mostly on the second question---the characterization of the stationary distributions of the iterative imputation conditional on its existence. Unlike usual MCMC algorithms, which are designed in such a way that the invariant distribution and target distribution are identical, the invariant distribution of iterative imputation (even if it exists) is largely unknown.

The analysis of iterative imputation is challenging for at least two reasons. First, the range of choices of conditional models is wide, and it would be difficult to provide a solution applicable to all situations. Second, the literature on Markov chains focuses on known transition distributions.  With iterative imputation, the distributions for the imputations are known only within specified parametric families.  For example, if a particular variable is to be updated conditional on all the others using logistic regression, the actual updating distribution depends on the logistic regression coefficients which are themselves estimated given the latest update of the missing values.

The main contribution of this paper is to develop a mathematical framework under which the asymptotic properties of iterative imputation can be discussed. In particular, we demonstrate the following results.
\begin{enumerate}
  \item \label{ContribComp}Given the existence of a unique invariant (stationary) distribution of the iterative imputation Markov chain, we provide a set of conditions under which this distribution converges in total variation to the posterior distribution of a joint Bayesian model, as the sample size tends to infinity. Under these conditions, iterative imputation is asymptotically equivalent to full Bayesian imputation using some joint model.  Among these conditions, the most important is that the conditional models are \emph{compatible}---that there exists a joint model whose conditional distributions are identical to the conditional models specified by the iterative imputation (Definition \ref{DefCompatible}). We discuss in Section \ref{SecCompatible}.

  \item \label{ContribCompNec}We consider model compatibility as a typically necessary condition for the iterative imputation distribution to converge to the posterior distribution of some Bayesian model (Section \ref{SecNec}).

  \item \label{ContribIncomp}For \emph{incompatible} models whose imputation distributions are generally different from any Bayesian model, we show that the combined completed-data maximum likelihood estimate of the iterative imputation is a consistent estimator if the set of conditional models is valid, that is, if each conditional family contains the true probability distribution (Definition \ref{DefValid} in Section \ref{SecIncomp}.).

\end{enumerate}

The analysis presented in this paper connects to the existing separate literatures on missing data imputation and Markov chain convergence.  Standard textbooks on imputation inference are \cite{LiRu02,Rubin1987}, and some key papers are \cite{Li1991,BaRu99,Meng1994,Meng1992,Rubin77,Rubin96,Sch97}. Large sample properties are studied by \cite{RoWa00,ScWe88,RoWa98}, small samples are by \cite{BaRu99}, and the issue of congeniality between the imputer's and analyst's models is considered by \cite{Meng1994}.

Our asymptotic findings for compatible and incompatible models use results on convergence of Markov chains, a subject on which there is a vast literature on stability and rate of convergence (\cite{GeGe84,Amit91,AmGr91}).  
In addition, empirical diagnostics of Markov chains have been suggested by many authors, for instance, \cite{rhat}. 
For the analysis of compatible models, we need to construct a bound for convergence rate using renewal theory \cite{Bax05,MT93,Ros95}, which  has the advantage of not assuming the existence of an invariant distribution, which is naturally yielded by the minorization and drift conditions.

In Section \ref{SecBack} of this article, we lay out our notation and assumptions. We then briefly review the framework of iterative imputation and the Gibbs sampler. In Section \ref{SecCompatible}, we investigate compatible conditional models. In Section \ref{SecIncomp}, the discussion focuses on incompatible models.  Section \ref{SecLinear} includes several simulation examples. An appendix is attached containing the technical developments and a brief review of the literature for Markov chain convergence via renewal theory.

\section{Background}\label{SecBack}
Consider a data set with $n$ cases and $p$ variables, where $\xx = (\xx_1,...,\xx_p)$ represents the complete data and $\xx_i = (x_{1,i},...,x_{n,i})^\top$
is the $i$-th variable. Let $\rr_i$ be the vector of observed data
indicators for variable $i$, equaling $1$ for observed variables and $0$ for missing, and let $\xx_i^{obs}$ and $\xx_i^{mis}$ denote
the observed and  missing subsets of variable $i$:
$$\xx^{obs}=\{\xx_i^{obs}: i=1,...,p\},\quad
\xx^{mis}=\{\xx_i^{mis}: i=1,...,p\}, \quad \rr = \{\rr_i: i =1,...,p\}.$$
To facilitate our
description of the procedures, we define
$$\xx^{obs}_{-j}=\{\xx_i^{obs}: i=1,...,j-1,j+1,...,p\},\quad
\xx^{mis}_{-j}=\{\xx_i^{mis}: i=1,...,j-1,j+1,...,p\}.$$ We use boldface $\xx$ to denote the entire data set and $x$ to denote individual observations. Therefore, $x_j$ denotes the $j$-th variable of one observation and $x_{-j}$ denotes all the variables except for the $j$-th one.

Throughout, we assume that the missing data process is ignorable. One set of sufficient conditions for ignorability is
that the $\rr_i$ process is missing at random and the parameter spaces for $\rr_i$ and $\xx$ are distinct, with independent
prior distributions \cite{LiRu02,Rubin1987}.

\subsection{Inference using multiple imputations}
Multiple imputation is a convenient tool to handle incomplete
data set by means of complete-data procedures. The framework
consists of producing $m$ copies of the imputed data and
applying the users' complete data procedures to each of
the multiply imputed data sets. Suppose that $m$ copies of point
estimates and variance estimates are obtained, denoted by
$(\hat\theta^{(i)},U^{(i)})$, $i=1,...,m$. The next step is to combine them into a single point estimate and a single variance estimate $(\hat \theta_m, \hat T_m)$ \cite{LiRu02}. If the imputed
data are drawn from the joint posterior distribution of the missing data under a
Bayesian model, under appropriate congeniality conditions,
$\hat \theta_m$ is asymptotically equal to the posterior mean of
$\theta$ and $\hat T_m$ is asymptotically equal to the posterior
variance of $\theta$ (\cite{Rubin1987,Meng1994}). The large
sample theory of Bayesian inference ensures that the posterior
mean and variance are asymptotically equivalent to the maximum
likelihood estimate and its variance based on the observed data
alone (see \cite{CoxHin74}). Therefore, the combined estimator
from imputed samples is efficient.
Imputations can also be constructed and used under other inferential frameworks; for example, Robins and Wang \cite{RoWa00,RoWa98} propose estimates based on estimating equations and derive
corresponding combining rules. For our purposes here, what is relevant is that the multiple imputations are being used to represent uncertainty about the joint distribution of missing values in a multivariate dataset.

\subsection{Bayesian modeling, imputation, and Gibbs sampling}\label{SecFull}
In Bayesian inference, multiply imputed data sets are treated as samples from the posterior distribution of the full (incompletely-observed) data matrix. In the parametric Bayesian approach, one specifies a family of distributions $f(\xx|\theta)$ and a prior $\pi(\theta)$ and then performs inference using  i.i.d.~samples from the
posterior predictive distribution,
\begin{equation}\label{Imp}p(\xx^{mis}|\xx^{obs}) = \int_{\Theta}
f(\xx^{mis}|\xx^{obs},\theta)p(\theta|\xx^{obs})d\theta,
\end{equation}
where $p(\theta|\xx)\propto \pi(\theta) f(\xx|\theta)$. Direct simulation from (\ref{Imp}) is generally difficult. One
standard solution is to draw approximate samples using the Gibbs sampler or some more complicated Markov chain Monte Carlo (MCMC) algorithm.  In
the scenario of missing data, one can use the ``data augmentation'' strategy to iteratively draw $\theta$ given
$(\xx^{obs},\xx^{mis})$ and $\xx^{mis}$ given $(\xx^{obs},\theta)$. Under regularity conditions (positive recurrence,
irreducibility, and aperiodicity; see \cite{GeGe84}), the Markov process is ergodic with limiting distribution $p(\xx^{mis},\theta|\xx^{obs})$.

In order to connect these results to the iterative imputation that is the subject of the present article, we consider a slightly different Gibbs
scheme which consists of indefinite iteration of following $p$ steps:,
\begin{description}
\item[Step 1.] Draw $\theta\sim p(\theta | \xx^{obs}_1,\xx_{-1})$ and
$\xx^{miss}_1\sim f(\xx^{miss}_1| \xx^{obs}_1,\xx_{-1},\theta)$;

\item[Step 2.] Draw $\theta \sim p(\theta|\xx^{obs}_2,\xx_{-2})$ and
$\xx^{miss}_2\sim f(\xx^{miss}_2| \xx^{obs}_2,\xx_{-2},\theta)$;

\quad \vdots

\item[Step $p$.] Draw $\theta \sim p(\theta| \xx^{obs}_p,\xx_{-p})$ and
$\xx^{miss}_p\sim f(\xx^{miss}_p| \xx^{obs}_p,\xx_{-p},\theta)$.
\end{description}
At each step, the posterior distribution is based on the updated values of the parameters and imputed data. It is not hard to verify that the Markov chain evolving according to steps 1 to $p$ (under mild regularity conditions) converges to the posterior distribution of the corresponding Bayesian model.

\subsection{Iterative imputation and compatibility}\label{SecIter} For  iterative imputation, we need to
specify $p$ conditional models,
$$g_j(\xx_j|\xx_{-j},\theta_j),$$
for $\theta_j \in \Theta_j$ with prior distributions $\pi_j (\theta_j)$ for $j=1,...,p$. Iterative imputation adopts the following scheme to construct a Markov chain,
\begin{description}
\item[Step 1.] Draw $\theta_1$ from $p_1 (\theta_1 | \xx^{obs}_1,\xx_{-1})$, which is
the posterior distribution associated with $g_1$ and $\pi_1$;
draw $\xx^{miss}_1$ from $g_1(\xx^{miss}_1|
\xx^{obs}_1,\xx_{-1},\theta_1)$;

\item[Step 2.] Draw $\theta_2$ from $p_2 (\theta_2 | \xx^{obs}_2,\xx_{-2})$, which is
the posterior distribution associated with $g_2$ and $\pi_2$;
draw $\xx^{miss}_2$ from $g_2(\xx^{miss}_2|
\xx^{obs}_2,\xx_{-2},\theta_2)$;

\quad \vdots

\item[Step $p$.] Draw $\theta_p$ from $p_p (\theta_p | \xx^{obs}_p,\xx_{-p})$, which is
the posterior distribution associated with $g_p$ and $\pi_p$;
draw $\xx^{miss}_p$ from $g_p(\xx^{miss}_p|
\xx^{obs}_p,\xx_{-p},\theta_p)$.
\end{description}

Iterative imputation has the practical advantage that, at each step, one only needs to set up a sensible regression model of $\xx_{j}$ given $\xx_{-j}$. This substantially reduces the modeling task, given that there are usually standard linear or generalized linear models for univariate responses of different variable types.
In contrast, full Bayesian (or likelihood) modeling requires the more difficult task of constructing a joint model for $\xx$.  Whether it is preferable to perform $p$ easy task or one difficult task, depends on the problem at hand.  All that is needed here is the recognition that, in {\em some} settings, users prefer the $p$ easy steps of iterative imputation.


But iterative imputation has conceptual problems.  Except in some special cases, there will not in general exist a
joint distribution of $\xx$ such that $f(\xx_j |
\xx_{-j},\theta) = g_j (\xx_j|\xx_{-j},\theta_j)$ for each $j$.
In addition, it is unclear whether the Markov process has a
probability invariant distribution; if there is such a distribution, it lacks characterization.

In this paper, we discuss the properties of the stationary distribution of the iterative imputation
process by first classifying the set of conditional models as compatible (defined as there existing a joint model $f$ which is consistent with all the conditional models) or incompatible.

We refer to the Markov chain generated by the scheme in Section
\ref{SecFull} as the \emph{Gibbs chain} and that generated by the
scheme in Section \ref{SecIter} as the \emph{iterative chain}. Our central analysis works by coupling the two.

\section{Compatible conditional models}\label{SecCompatible}

\subsection{Model compatibility}
Analysis of iterative imputation is particularly challenging
partly because of the large collection of possible choices of
conditional models.
We begin by considering a restricted class, \emph{compatible conditional
models}, defined as follows:

\begin{definition}\label{DefCompatible}
A set of conditional models $\{g_j(x_j|x_{-j},\theta_j):
\theta_j \in \Theta_j, j = 1,...,p\}$ is said to be
\emph{compatible} if there exists a joint model $\{f(x|\theta): \theta\in \Theta\}$
and a collection of surjective maps, $\{t_j : \Theta \rightarrow \Theta_j
: j=1,...,p\}$ such that for each $j$, $\theta_j \in \Theta_j$,
and $\theta \in t_j^{-1}(\theta_j)= \{\theta: t_j(\theta) = \theta_j\}$,
$$g_j(x_j | x_{-j},\theta_j)=f(x_j | x_{-j},\theta).$$
Otherwise, $\{g_j: j =1,...,p\}$ is said to be \emph{incompatible}.
\end{definition}
Though imposing certain restrictions, compatible models do include quite a collection of procedures practically in use (e.g.\ {\tt ice} in \texttt{Stata}). In what
follows, we give a few examples of compatible and
incompatible conditional models.

We begin with a simple linear model, which we shall revisit in Section \ref{SecLinear}.

\begin{example}[bivariate Gaussian]\label{ExToy}
Consider a binary continuous variable $(x,y)$ and conditional models
$$x|y \sim N(\alpha_{x|y}+\beta_{x|y}y, \tau^2_x),\quad y|x \sim N(\alpha_{y|x}+\beta_{x|y}x, \tau^2_y).$$
These two conditional models are compatible if and only if $(\beta_{x|y},\beta_{x|y},\tau_x,\tau_y)$ lie on a subspace determined from the joint model,
\[
\left(
\begin{array}{c}
x \\
y%
\end{array}%
\right) \sim N\left( \left(
\begin{array}{c}
\mu_{x} \\
\mu_{y}%
\end{array}%
\right) ,\Sigma \right) ,\quad
\mbox{where  }
\Sigma =\left(
\begin{array}{cc}
\sigma_{x}^{2} & \rho \sigma_{x}\sigma_{y} \\
\rho \sigma_{x}\sigma_{y} & \sigma_{y}^{2}%
\end{array}%
\right),
\]
with $\sigma_x, \sigma_y >0$ and $\rho \in [-1,1]$.
The reparameterization from $(\mu_x,\mu_y,\sigma_x,\sigma_y,\rho)$ to the parameters of the conditional models is:
\begin{eqnarray*}t_{1}(\mu_{x},\sigma_{x}^{2},\mu_{y},\sigma_{y}^{2},\rho
)=(\alpha_{x|y},\beta_{x|y},\tau_x^2)=\Big(\mu_{x}-\frac{\rho \sigma
_{x}}{\sigma_{y}}\mu_{y},\frac{\rho \sigma_{x}}{\sigma_{y}},(1-\rho
^{2})\sigma_{x}^{2}\Big)\\
t_{2}(\mu_{x},\sigma_{x}^{2},\mu_{y},\sigma_{y}^{2},\rho
)=(\alpha_{y|x},\beta_{y|x},\tau_y^2)=\Big(\mu_{y}-\frac{\rho \sigma
_{y}}{\sigma_{x}}\mu_{x},\frac{\rho \sigma_{y}}{\sigma_{x}},(1-\rho
^{2})\sigma_{y}^{2}\Big ).
\end{eqnarray*}

\end{example}

\noindent
The following example is a natural extension.

\begin{example}[continuous data]\label{ExLinear}
Consider a set of conditional linear models:  for each $j$,
$$x_j | x_{-j}, \beta_j, \sigma^2_j \sim N\left((\mathbf
1,x_{-j})\beta_j,\sigma^2_j \right),$$ where $\beta_j$ is a
$p\times 1$ vector, $\mathbf 1 = (1,...,1)^\top$. Consider
the joint model of $(x_1,...,x_p)\overset {i.i.d.}\sim N(\mu,
\Sigma)$. Then the conditional distribution of each $x_j$ given
$x_{-j}$ is Gaussian. The maps $t_j$'s can be derived by
conditional multivariate Gaussian calculations.
\end{example}


\begin{example}[continuous and binary data]\label{ExLogit}
Let $x_1$ be a Bernoulli random variable and $x_2$ be a
continuous random variable. The conditional models are as
follows:
$$x_1|x_2 \sim Bernoulli \left(\frac{e^{\alpha + \beta
x_2}}{1+e^{\alpha + \beta x_2}}\right),\quad x_2|x_1
\sim N(\beta_0 + \beta_1 x_1 ,\sigma^2).$$ The above conditional
models are compatible with the following joint model:
$$x_1 \sim Bernoulli (p), \quad x_2|x_1
\sim N(\beta_0 + \beta_1 x_1 ,\sigma^2).$$
If we let
\begin{eqnarray*}
t_1 (p,\beta_0,\beta_1,\sigma^2)&=&\left(\log\frac p {1-p}
-\frac{\beta_1^2}{2\sigma^2}, \frac {\beta_1}{2
\sigma^2}\right)=(\alpha,\beta)\\
t_2 (p,\beta_0,\beta_1,\sigma^2)&=& (\beta_0,\beta_1),
\end{eqnarray*}
the conditional models and this joint model are
compatible with each other. Similarly compatible models can be
defined for other natural exponential families. See \cite{Efron75,MN98}.
\end{example}

\begin{example}[incompatible Gaussian conditionals]\label{ExCounter}
There are many incompatible conditional models. For
instance, $$x|y \sim N(\beta_1 y+ \beta_2 y^2, 1),\quad y|x
\sim N(\lambda_1 x, 1),$$
are compatible only if $\beta_2 = 0$.
\end{example}

%
%

\subsection{Total variation distance between  two transition kernels}\label{SecPriorImp}

Let $\{\xx^{mis,1}(k): k \in \mathbb Z^+\}$ be the Gibbs chain and $\{\xx^{mis,2}(k): k \in \mathbb Z^+\}$ be the iterative chain.
Both chains live on the space of the missing data. 
We write the completed data as $\xx^{i}(k) = (\xx^{mis,i}(k),\xx^{obs})$ for the Gibbs chain ($i=1$) and the iterative chain ($i=2$).  The transition kernels are
\begin{equation}  \label{kernel}
K_i (w,dw')=P(\xx^{mis,i}(k+1) \in d w'|\xx^{mis,i}(k)=w), \mbox{ for } i=1,2.
\end{equation}
where $w$ is a generic notation for the state of the processes.
The transition kernels ($K_1$ and $K_2$) depend on $\xx^{obs}$. For simplicity, we
omit the index of $\xx^{obs}$ in the notation of $K_i$. Also, we let
$$K_i^{(k)} (\nu, A)\triangleq P_\nu(\xx^{mis,i}(k)\in A),$$
for $\xx^{mis,i}(0)\sim \nu$, $\nu$ being some starting distribution.
The probability measure $P_{\nu}$ also depends on $\xx^{obs}$.
Let $d_{TV}$ denote the total variation distance between two measures, that is, for two measures, $\nu_1$ and $\nu_2$, defined on the same probability space
$$d_{TV}(\nu_1,\nu_2) = \sup_{A\in \mathcal F} |\nu_1 (A) - \nu_2 (A)|.$$
We further define
$$\Vert\nu\Vert_V = \sup_{|h|\leq V}\int h(x) \nu (dx)$$
and $\Vert\nu\Vert_{\mathbf 1} = \Vert\nu\Vert_{V} $ for $V\equiv 1$.
Let $\nu^{\xx^{obs}}_{i}$ be the stationary distribution of $K_{i}$. We intend to establish conditions under which 
$$d_{TV}(\nu^{\xx^{obs}}_{1}, \nu^{\xx^{obs}}_{2})\rightarrow 0$$
in probability as $n\rightarrow \infty$ and thus the iterative imputation and the joint Bayesian imputation are asymptotically the same.


Our basic strategy for analyzing the compatible conditional models is to first establish that the transition kernels $K_{1}$ and $K_{2}$ are close to each other in a large region $A_{n}$ (depending on the observed data $\xx^{obs}$), that is, $\Vert K_{1}(w,\cdot) - K_{2}(w,\cdot)\Vert_{\mathbf 1}\rightarrow 0$ as $n\rightarrow \infty$ for $w\in A_{n}$; and, second, to show that the two stationary distributions are close to each other in total variation in that the stationary distributions are completely determined by the transition kernels. In this subsection, we start with the first step, that is, to show that $K_{1}$ converges to $K_{2}$.


Both the Gibbs chain and the iterative chain evolve by updating each missing variable from the corresponding posterior predictive distributions. Upon comparing the difference between
the two transition kernels associated with the simulation schemes in Sections \ref{SecFull} and \ref{SecIter}, it suffices to compare the following  posterior predictive distributions (for each $j=1,...,p$),
\begin{eqnarray}
f(\xx_{j}^{mis}|\xx^{obs}_j,\xx_{-j}) &=&\int f(\xx_{j}^{mis}|\xx^{obs}_j,\xx_{-j},\theta )p(\theta |\xx^{obs}_j,\xx_{-j})d\theta  \label{f} \\
g_{j}(\xx_{j}^{mis}|\xx^{obs}_j,\xx_{-j}) &=&\int g_{j}(\xx_{j}^{mis}|\xx^{obs}_j,\xx_{-j},\theta_{j})p_{j}(\theta_{j}|
\xx^{obs}_j,\xx_{-j})d\theta_{j},
\label{g}
\end{eqnarray}%
where $p$ and $p_j$ denote the posterior distributions under $f$ and $g_j$ respectively.
Due to compatibility, the distributions of the missing data given the parameters are the same for
the joint Bayesian model and the iterative imputation model:
\begin{equation*}
f(\xx_{j}^{mis}|\xx^{obs}_j,\xx_{-j},\theta )=g_{j}(\xx_{j}^{mis}|\xx^{obs}_j,\xx_{-j},\theta_{j}),
\end{equation*}
if $t_{j}(\theta )=\theta_{j}$. The only difference lies in their posterior
distributions. In fact,  the $\Vert \cdot \Vert _{\mathbf 1}$ distance between two posterior predictive distributions is bounded by the distance between the posterior distributions of parameters. Therefore, we move to comparing $p(\theta|\xx^{obs}_j,\xx_{-j})$ and $p_j(\theta_j|\xx^{obs}_j,\xx_{-j})$.

\paragraph{Parameter augmentation.}

Upon comparing the posterior distributions of $\theta$ and $\theta_{j}$, the
first disparity to reconcile is that the dimensions are usually different.
Typically $\theta_{j}$ is of a lower dimension. Consider the linear model in
Example \ref{ExToy}. The conditional models include three parameters (two
regression coefficients and variance of the errors), while the joint model
has five parameters $(\mu_{x},\mu_{y},\sigma_{x},\sigma_{y},\rho )$.
This is because the (conditional) regression models are usually conditional
on the covariates. The joint model not only parameterizes the conditional
distributions of $\xx_{j}$ given $\xx_{-j}$
but also the marginal distribution of $\xx_{-j}$. Therefore, it includes extra
parameters, although the distributions of the missing data is independent of
these parameters. We augment the parameter space of the iterative imputation
to $(\theta_{j},\theta_{j}^{\ast })$ with the corresponding map $\theta
_{j}^{\ast }=t_{j}^{\ast }(\theta )$. The augmented parameter $(\theta_{j},\theta_{j}^{\ast })$ is
a non-degenerated reparameterization of $\theta $, that is, $T_{j}(\theta
)=(t_{j}(\theta ),t_{j}^{\ast }(\theta ))$ is a one-to-one (invertible) map.

To illustrate this parameter augmentation, we consider the linear model in
Example \ref{ExToy} in which $\theta =(\mu_{x},\sigma_{x}^{2},\mu
_{y},\sigma_{y}^{2},\rho )$, where we use $\mu_{x}$ and $\sigma_{x}^{2}$
to denote mean and variance of the first variable, $\mu_{y}$ and $\sigma
_{y}^{2}$ to denote the mean and variance of the second, and $\rho $ to
denote the correlation. The
reparameterization is,
\begin{eqnarray*}
\theta_{2} &=&t_{2}(\mu_{x},\sigma_{x}^{2},\mu_{y},\sigma_{y}^{2},\rho
)=(\alpha_{y|x},\beta_{y|x},\tau_y^2)=(\mu_{y}-\frac{\rho \sigma
_{y}}{\sigma_{x}}\mu_{x},\frac{\rho \sigma_{y}}{\sigma_{x}},(1-\rho
^{2})\sigma_{y}^{2}),\\
\theta_{2}^{\ast } &=&t_{2}^{\ast }(\mu_{x},\sigma_{x}^{2},\mu
_{y},\sigma_{y}^{2},\rho )=(\mu_{x},\sigma_{x}^{2}).
\end{eqnarray*}%
The function $t_{2}$ maps to the regression coefficients and the variance of the
residuals; $t_{2}^{\ast }$ maps to the marginal mean and variance of $x$. Similarly, we can define the map of $t_{1}$ and $t_{1}^{\ast
}$.

\paragraph{Impact of the prior distribution.}

Because we are assuming compatibility, we can drop the notation $g_{j}$ for conditional model of the $j$-th variable. Instead, we unify the
notation to that of the joint Bayesian model $f(\xx_{j}|%
\xx_{-j},\theta )$. In addition, we abuse the notation and write $f(%
\xx_{j}|\xx_{-j},\theta_{j})=f(%
\xx_{j}|\xx_{-j},\theta )$ for $t_{j}(\theta )=\theta_{j}$.
For instance, in Example \ref{ExToy}, 
we write $f(y|x,\alpha_{y|x},\beta_{y|x},\sigma_{y|x})=f(y|x,\mu_{x},\mu_{y},\sigma_{x},\sigma_{y},\rho )$
as long as $\alpha_{y|x}=\mu_{y}-\frac{\rho \sigma_{y}}{\sigma_{x}}\mu
_{x}$, $\beta_{y|x}=\frac{\rho \sigma_{y}}{\sigma_{x}}$, and $\sigma
_{y|x}^{2}=(1-\rho ^{2})\sigma_{y}^{2}$.

The prior distribution $\pi$ on $\theta$ for the joint Bayesian model implies a prior
on $(\theta_{j},\theta_{j}^{\ast })$, denoted by
\begin{equation*}
\pi_{j}^{\ast }(\theta_{j},\theta_{j}^{\ast })=|\det (\partial
T_{j}/\partial \theta )|^{-1}\pi (T_{j}^{-1}(\theta_{j},\theta_{j}^{\ast
})).
\end{equation*}%
For the full Bayesian model, the posterior distribution of $\theta_{j}$ is
\begin{equation*}
p(\theta_{j}|\xx^{obs}_j,\xx_{-j})=\int
p(\theta_{j},\theta_{j}^{\ast }|\xx^{obs}_j,\xx_{-j})d\theta_{j}^{\ast }\propto \int f(\xx^{obs}_j,\xx_{-j}|\theta_{j},\theta_{j}^{\ast })\pi
_{j}^{\ast }(\theta_{j},\theta_{j}^{\ast })d\theta_{j}^{\ast }.
\end{equation*}%
Because $f(\xx_{j}^{obs}|\xx_{-j},\theta
_{j},\theta_{j}^{\ast })=f(\xx_{j}^{obs}|\xx%
_{-j},\theta_{j})$, the above posterior distribution can be further
reduced to
\begin{equation*}
p(\theta_{j}|\xx^{obs}_j,\xx_{-j})%
\propto f(\xx_{j}^{obs}|\xx_{-j},\theta_{j})\int f(\xx_{-j}|\theta_{j},\theta
_{j}^{\ast })\pi_{j}^{\ast }(\theta_{j},\theta_{j}^{\ast })d\theta
_{j}^{\ast }.
\end{equation*}%
If we write
\begin{equation*}
\pi_{j,\xx_{-j}}(\theta_{j})\triangleq  \int f(%
\xx_{-j}|\theta_{j},\theta_{j}^{\ast })\pi_{j}^{\ast }(\theta_{j},\theta
_{j}^{\ast })d\theta_{j}^{\ast },
\end{equation*}%
then the posterior distribution of $\theta_j$ under the joint Bayesian model is
\begin{equation*}
p(\theta_{j}|\xx^{obs}_j,\xx_{-j})%
\propto f(\xx_{j}^{obs}|\xx_{-j},\theta_{j})\pi_{j,\xx_{-j}}(\theta_{j}).
\end{equation*}%
Compared with the posterior distribution of the iterative imputation
procedure, which is proportional to
\begin{equation*}
p_{j}(\theta_{j}|\xx^{obs}_j,\xx_{-j})\propto g_{j}(\xx_{j}^{obs}|\xx%
_{-j},\theta_{j})\pi_{j}(\theta_{j})=f(\xx_{j}^{obs}|%
\xx_{-j},\theta_{j})\pi_{j}(\theta_{j}),
\end{equation*}%
the difference lies in the prior distributions, $\pi_{j}(\theta_{j})$ and $%
\pi_{j,\xx_{-j}}(\theta_{j})$.

\paragraph{Controlling the  distance between the posterior
predictive distributions.}

We put forward tools to control the  distance between the two posterior predictive distributions in \eqref{f} and \eqref{g}.  Let $\xx$ be the generic notation the observed data, and let $f_{\xx}(\theta)$ and $g_{\xx}(\theta)$ be two posterior densities of $\theta$. Let $h(\tilde{x}|\theta )$ be the density function for
future observations given the parameter $\theta $, and let $\tilde{f}_{\xx}(\tilde{x})$ and $\tilde{g}_{\xx}(\tilde{x})$ be the posterior predictive distributions:
\begin{equation*}
\tilde{f}_{\xx}(\tilde{x})=\int h(\tilde{x}|\theta )f_{%
\xx}(\theta )d\theta ,\quad \tilde{g}_{\xx}(%
\tilde{x})=\int h(\tilde{x}|\theta )g_{\xx}(\theta )d\theta .
\end{equation*}%
It is straightforward to obtain that
\begin{equation}\label{PropPred}
\Vert\tilde{f}_{\xx}-\tilde{g}_{\xx}\Vert_{\mathbf 1}\leq
\Vert f_{\xx}-g_{\xx}\Vert_{\mathbf 1}.
\end{equation}
The next proposition provides sufficient conditions that $\Vert f_{\xx}-g_{\xx}\Vert_{\mathbf 1}$ vanishes.

\begin{proposition}
\label{PropPrior} Let $n$ be the sample size. Let $f_{\xx%
}(\theta )$ and $g_{\xx}(\theta )$ be two posterior density
functions that share the same likelihood but have two different prior
distributions $\pi_{f}$ and $\pi_{g}$. Let
$$L(\theta )=\frac{\pi_{g}(\theta)}{\pi_{f}(\theta )},\qquad r(\theta) = \frac{g_{\xx}(\theta)}{f_{\xx}(\theta)} = \frac{L(\theta)}{\int L(\theta) f_{\xx}(\theta)d \theta},
$$ and $n$ denote sample size.
Let $\partial L(\theta)$ be the partial derivative with respect to $\theta$ and let $\xi $ be a random variable such that
\begin{equation*}
L(\theta )=L(\mu_{\theta })+\partial L(\xi )\cdot (\theta -\mu_{\theta }),
\end{equation*}%
where `` $\cdot$'' denotes inner product and $\mu_{\theta} = \int \theta f_{\xx}(\theta)d\theta$.
If there exists a random variable $Z(\theta)$ with finite variance under $f_{\xx}$, such that 
\begin{equation}
\left\vert \sqrt n \partial L(\xi )\cdot (\theta -\mu
_{\theta})\right\vert\leq |\partial L(\mu_{\theta})|Z(\theta),  \label{Dom}
\end{equation}
then there exists a constant $\kappa>0$ such that for $n$ sufficiently large
%
\begin{equation}\label{pbd}\Vert\tilde{f}_{\xx}-\tilde{g}_{\xx}\Vert_{\mathbf 1}\leq \frac{\kappa \sqrt {|\partial \log L(\mu_{\theta})|}}{n^{1/4}}.
\end{equation}
\end{proposition}
We prove this proposition in Appendix \ref{Apdbd}.

\begin{remark}\label{RemLocal}
We adapt Proposition \ref{PropPrior} to the analysis of the conditional models. Expresion \eqref{Dom}  implies that the posterior standard deviation of $\theta$ is $O(n^{-1/2})$.
For most parametric models, \eqref{Dom} is satisfied as long as the observed Fisher information is  bounded from below by  $\varepsilon n$ for some $\varepsilon >0$.
In particular, we let $\hat \theta(\xx)$ be the complete-data MLE and $A_{n} = \{\xx: |\hat \theta(\xx)| \leq \gamma\}$.  Then, \eqref{Dom} is satisfied on the set $A_{n}$ for any fixed $\gamma$.
\end{remark}

\begin{remark}
In order to verify that $\partial\log L(\theta)$ is bounded, one only needs to know $\pi_f$ and $\pi_g$ up to a normalizing constant. This is because the bound is in terms of $\partial L(\theta)/L(\theta)$. This helps to handle the situation when improper priors are used and it is not feasible to obtain a normalized prior distribution.
In the current context, the prior likelihood ratio is $L(\theta_{j})= \pi_{j}(\theta_{j})/\pi_{j,\xx_{-j}}(\theta_{j})$.
\end{remark}

\begin{remark}
If $L(\theta)$ is twice differentiable, the convergence rate in \eqref{pbd} can be improved to $O(n^{-1/2})$. However, $O(n^{-1/4})$ is sufficient for the current analysis.
\end{remark}



\subsection{Convergence of the invariant distributions} \label{SecCgt}

With Proposition \ref{PropPrior} and Remark \ref{RemLocal}, we have established that the transition kernels of the Gibbs chain and the iterative chain are close to each other in a large region $A_{n}$.
The subsequent analysis falls into several steps. First, we slightly modify the processes by conditioning them on the set $A_{n}$ with stationary distributions $\tilde \nu^{\xx^{obs}}_{i}$ (details provided below).
The  stationary distributions of the conditional processes and the original processes ($\tilde \nu^{\xx^{obs}}_{i}$ and $\nu^{\xx^{obs}}_{i}$) are close in total variation. Second, we show (in Lemma \ref{ThmComp}) that, with a bound on the  convergence rate, $\tilde \nu^{\xx^{obs}}_{1}$ and $\tilde \nu^{\xx^{obs}}_{2}$ are close in total variation and so it is with $\nu^{\xx^{obs}}_{1}$ and $ \nu^{\xx^{obs}}_{2}$. The bound of convergence rate can be established  by Proposition \ref{PropDrift}.

To proceed, we consider the chains conditional on the set $A_{n}$ where the two transition kernels are uniformly  close to each other. In particular, for each set $B$, we let
\begin{equation}\label{condK}
\tilde K_{i} (w,B) = \frac { K_{i}(w,B\cap A_{n})}{ K_{i}(w,A_{n})}.
\end{equation}
That is, we create another two processes, for which we update the missing data conditional on $\xx \in A_{n}$.
The next lemma shows that we only need to consider the chains conditional on the set $A_{n}$.

\begin{lemma}\label{LemCond}
Suppose that both $K_{1}$ and $K_{2}$ are positive Harris recurrent.
We can choose $A_{n}$ as in the form of Remark \ref{RemLocal} and $\gamma$ sufficiently large so that
\begin{equation}\label{condcgt}
\nu^{\xx^{obs}}_{i}(A_{n})\rightarrow 1
\end{equation}
in probability as $n\rightarrow \infty$.
Let $\tilde {\xx}^{mis,i}(k)$ be the Markov chains following $\tilde K_{i}$, defined as in \eqref{condK}, with invariant distribution $\tilde\nu^{\xx^{obs}}_{i}$. Then,
\begin{equation}
d_{TV}(\nu^{\xx^{obs}}_{i},\tilde \nu^{\xx^{obs}}_{i})\rightarrow 0,
\end{equation}
as $n\rightarrow \infty$.
\end{lemma}

The proof is elementary by the representation of $\nu_{i}^{\xx^{obs}}$ through the renewal theory and therefore is omitted.
Based on the above lemma, we only need to show that $d_{TV}( \tilde \nu_{1}^{\xx^{obs}},\tilde \nu^{\xx^{obs}}_{2})\rightarrow 0$.   
The expression $\Vert K_{1}(w,\cdot)-K_{2}(w,\cdot)\Vert_{\mathbf 1}$ approaches 0 uniformly for $w\in A_{n}$.
This implies that 
\begin{equation*}
\Vert\tilde K_{1}(w,\cdot),\tilde K_{2}(w,\cdot)\Vert_{\mathbf 1}\rightarrow 0
\end{equation*}
as $n\rightarrow\infty$ uniformly for $w\in A_{n}$.
With the above convergence, we use the following lemma to establish the convergence between $\tilde \nu_{1}^{\xx^{obs}}$ and $\tilde \nu^{\xx^{obs}}_{2}$.

\begin{lemma}
\label{ThmComp} Let $\tilde{\xx}^{mis,i}(k)$ admit data-dependent transition kernels $\tilde K_i$ for $i=1,2$. We use $n$ to denote sample size. Suppose that each $\tilde K_i$ admits a data-dependent  unique invariant distribution, denoted by $\tilde \nu^{\xx^{obs}}_i$, and that the following two conditions hold:
\begin{enumerate}
\item The convergence of the two transition kernels
\begin{equation}
d(A_{n})\triangleq\sup_{w\in A_{n}}\Vert\tilde K_{1}(w,\cdot )-\tilde K_{2}(w,\cdot )\Vert_{V}\rightarrow
0,  \label{BdKernal}
\end{equation}%
in probability as $n\rightarrow \infty $. The function  $V$ is either a geometric drift function for $ \tilde K_2$ or a constant, i.e., $V=1$.
\item
Furthermore, there exists a monotone decreasing sequence $r_k \rightarrow 0$ (independent of data) and a starting measure $\nu$ (depending on data) such that
\begin{equation}  \label{BdProb}
P\left[\Vert \tilde K_i^{(k)}(\nu,\cdot)-\tilde \nu^{\xx^{obs}}_i(\cdot)\Vert_V\leq r_k,\forall
k>0\right] \rightarrow 1,
\end{equation}
as $n\rightarrow \infty$.

\end{enumerate}
Then,
\begin{equation}  \label{CgtMea}
\Vert\tilde \nu^{\xx^{obs}}_1-\tilde \nu^{\xx^{obs}}_2\Vert_V\rightarrow 0,
\end{equation}
in probability as $n\rightarrow \infty$.
\end{lemma}

\begin{remark}
The above lemma holds if $V=1$ or $V$ is a drift function. For the analysis of convergence in total variation, we only need that $V=1$. The results when $V$ is a drift function is prepared for the analysis of incompatible models. 
\end{remark}

The first condition in the above lemma has been obtained by the result of Proposition \ref{PropPrior}.
Condition \eqref{BdProb} is more difficult to establish.
According to the standard results in \cite{Ros95} (see also \eqref{cgtbd} in the appendix),
one set of sufficient conditions for  \eqref{BdProb} is that the chains $\tilde K_{1}$ and $\tilde K_{2}$ admit a common  small set, $C$; in addition, each of them admits their own drift functions associated with the small set $C$ (c.f. Appendix \ref{SecMC}). 

Gibbs chains typically admit a small set $C$ and a drift function $V$, that is, for some positive measure $\mu$
\begin{equation}\label{small}
\tilde K_1(w,A)\geq q_{1} \mu_{1}(A),
\end{equation}
for $w\in C$, $q_{1}\in (0,1)$; for some $\lambda_{1}\in (0,1)$ and  for all $w\notin C$
\begin{equation}\label{drift1}
\lambda_1 V(w)\geq\int V(w') \tilde K_{1}(w,dw').
\end{equation}
With the existence of $C$ and $V$ a bound of convergence $r_{k}$ (with starting point $w\in C$) can be established for the Gibbs chain by standard results (see, for instance, \cite{Ros95}), and $r_{k}$ only depends on $\lambda_{1}$ and $q_{1}$. Therefore, it is necessary to require that $\lambda_{1}$ and $q_{1}$ are independent of $\xx^{obs}$. In contrast, the small set $C$ and drift function $V$ could be data-dependent.

Given that $\tilde K_{1}$ and $\tilde K_{2}$ are close in ``$\Vert \cdot \Vert_{\mathbf 1}$'', the set $C$ is also a  small set for $\tilde K_{2}$, that is $\tilde K_2(w,A)\geq q_{2} \mu_{2}(A),$ for some $q_{2}\in (0,1)$, all $w\in C$, and all measurable set $A$.
The following proposition, whose proof is given in the appendix, establishes the conditions under which  $V$ is also a drift function for $\tilde K_{2}$.



\begin{proposition}\label{PropDrift}
Assume the following conditions hold.
\begin{enumerate}
\item The transition kernel $\tilde K_{1}$ admits a small set $C$ and a drift function $V$  satisfying \eqref{drift1}.  
\item Let $L_{j}(\theta_{j}) = \pi_j(\theta_j)/\pi_{j,\xx_{-j}}(\theta_j)$ ($j=1,...,p$) be the the ratio of prior distributions for each conditional model (possibly depending on the data) so that on the set $A_{n}$
$\sup_{|\theta_{j}|<\gamma}\partial L_{j}(\theta_{j})/L_{j}(\theta_{j}) <\infty.$
\item
For each $j$ and $1\leq k \leq p-j$, there exists a $Z_{j}(\theta_{j})$ serving as the bound in \eqref{Dom} for each $L_{j}$. In addition, $Z_{j}$ satisfies the following moment condition
\begin{equation}\label{cond2}
 \tilde E_{1}\left[~ Z_{j+1}^{2}(\theta_{j+1}) V^{2}(w_{j+k})~|~w_{j} ~\right]= o( n)V^{2}(w_{j}), \qquad 
\end{equation}
where $ \tilde E_{1}$ is the expectation associated with the updating distribution of  $\tilde K_{1}$ and $w_{j}$ is the state of the chain when the $j$-th variable is just updated. The convergence $o(n)/n \rightarrow 0$ is uniform in $w_{j}\in A_{n}$.
\end{enumerate}
Then, there exits $\lambda_2\in (0,1)$ such that as $n$ tends to infinity with probability converging to one the following inequality holds
\begin{equation}\label{Lyp}\lambda_2 V(w)\geq\int V(w') \tilde K_{2}(w,dw').
\end{equation}

\begin{remark}
The intuition of the above proposition is as follows.
$V$ satisfying inequality \eqref{drift1} is a drift function of $\tilde K_{1}$ to $C$.
Since the $\tilde K_{1}$ and $\tilde K_{2}$ are close to each other, we may expect that 
$\int V(w') \tilde K_{1}(w,dw') \approx \int V(w') \tilde K_{2}(w,dw')$.
The above proposition basically states the conditions under which this approximation is indeed true and suggests that $V$ be a drift function of $\tilde K_{2}$ if it is a drift function of $\tilde K_{1}$.
Condition \eqref{cond2} is imposed for a technical purpose.
In particular, we allow the expectation of $ Z^{2}_{j+1}(\theta_{j+1}) V^{2}(w_{j+k})$ to grow to infinity but at a slower rate than $ n$. Therefore, it is a mild condition.
\end{remark}

%
%
%
\end{proposition}

We now summarize the analysis and the results of the compatible conditional models in the following theorem.

\begin{theorem}\label{ThmCompRev}
Suppose that a set of conditional models $\{g_{j}(x_{j}|x_{-j},\theta_{j}): \theta_{j}\in \Theta_{j},j=1,...,p\}$ is compatible with a joint model $\{f(x|\theta): \theta \in \Theta\}$.  The Gibbs chain and the iterative chain then admit transition kernels $K_{i}$ and unique stationary distributions $\nu^{\xx^{obs}}_{i}$. Suppose the following conditions are satisfied:
\begin{itemize}
\item [A1] Let $A_{n} = \{\xx: |\hat \theta(\xx)| \leq \gamma\}$. One can choose $\gamma$ sufficiently large so that
$\nu^{\xx^{obs}}_{i}(A_{n})\rightarrow 0,$
in probability as $n\rightarrow \infty$.

\item[A2] The conditions in Proposition \ref{PropDrift} hold.



\end{itemize}
 Then,
\begin{equation*}
d_{TV}(\nu_{1}^{\xx^{obs}}, \nu_{2}^{\xx^{obs}})\rightarrow 0
\end{equation*}
in probability as $n\rightarrow \infty$.
\end{theorem}

\begin{remark}
One sufficient condition for A1 is that the stationary distributions of $\hat \theta(\xx)$ under $\nu^{\xx^{obs}}_{i}$ converge to a value $\theta^{i}$, where $\theta^{1}$ and $\theta^{2}$ are not necessarily the same.
\end{remark}

\begin{remark}
In addition to the conditions of Proposition \ref{PropPrior}, Proposition \ref{ThmComp} 
also requires that one constructs a drift function towards a small set for the Gibbs chain. One can usually construct $q_{1}$ and $\lambda_{1}$ free of data if the proportion of missing data is  bounded from the above by  $1-\varepsilon$. 
The most difficult task usually lies in constructing a drift function. For illustration purpose, we construct a drift function (in the supplement material) for the linear example in Section \ref{SecLinear}.
\end{remark}

\begin{proof}[Proof of Theorem \ref{ThmCompRev}]
We summarize the analysis of compatible models in this proof. If $g_{j}$'s are compatible with $f$, then the  conditional posterior predictive distributions of the Gibbs chain and the iterative chain are given in \eqref{f} and \eqref{g}. Thanks to compatibility, the ``$\Vert \cdot\Vert_{\mathbf 1}$'' distance between the posterior predictive distributions are bounded by the distance between the posterior distributions of their own parameters as  in \eqref{PropPred}.

 On the set $A_{n}$, the Fisher information of the likelihood has a lower bound of $\varepsilon n$ for some $\varepsilon$. Then, by Proposition \ref{PropPrior} and the second condition in Proposition \ref{PropDrift}, the  distance between the two posterior distributions is of order $O(n^{-1/4})$ uniformly on set $A_{n}$. Similar convergence result holds for the conditional transition kernels, that is, $\Vert\tilde K_{1}(w,\cdot)-\tilde K_{2}(w,\cdot)\Vert_{\mathbf 1}\rightarrow 0$. 
Thus, the  first condition in Lemma \ref{ThmComp} has been satisfied. 
 

To verify the conditions of Proposition \ref{PropDrift}, one needs to construct a small set $C$ such that \eqref{small} holds for both chains and a drift function $V$ for one of the two chains such that \eqref{drift1} holds.
Based on  the results of Proposition \ref{PropDrift}, $\tilde K_{1}$ and $\tilde K_{2}$ share a common data-dependent small set $C$ with $q_{i}$ independent of data
and a drift function $V$ (possibly with different rate $\lambda_1$ and $\lambda_2$).

According to the standard bound of Markov chain rate of convergence (for instance, \cite{Ros95} and \eqref{cgtbd} in the appendix), there exists a common starting value $w\in C$ and a bound $r_{k}$ such that the bound \eqref{BdProb} in Lemma \ref{ThmComp} is satisfied.
Thus both conditions in  Lemma \ref{ThmComp} have been satisfied and  further
$$d_{TV}(\tilde\nu_{1}^{\xx^{obs}},\tilde\nu_{2}^{\xx^{obs}})\rightarrow 0,$$
in probability as $n\rightarrow \infty$.
According to condition A1 and Lemma \ref{LemCond}, the above convergence implies that
$$d_{TV}(\nu_{1}^{\xx^{obs}},\nu_{2}^{\xx^{obs}})\rightarrow 0.$$
Thereby, we conclude the analysis.
\end{proof}

\subsection{On the necessity of model
compatibility}\label{SecNec}

Theorem \ref{ThmCompRev} shows that for compatible models and under suitable technical conditions, iterative imputation is asymptotically equivalent to Bayesian imputation. The following theorem suggests that model compatibility is typically necessary for this convergence.

Let $P^{f}$ denote the probability measure induced by the posterior predictive distribution of the joint Bayesian model and $P^{g}_j$ denote those induced by the iterative imputation's conditional models. That is,
\begin{eqnarray*}P^{f}(\xx_j^{mis}\in A|\xx^{mis}_{-j},\xx^{obs})&=& \int_A f(\xx_{j}^{mis}|\xx^{mis}_{-j},\xx^{obs},\theta )p(\theta |\xx^{mis}_{-j},\xx^{obs})d\theta\\
P^{g}_j(\xx_j^{mis}\in A|\xx^{mis}_{-j},\xx^{obs})&=& \int_A g_j(\xx_{j}^{mis}|\xx^{mis}_{-j},\xx^{obs},\theta )p_j(\theta |\xx^{mis}_{-j},\xx^{obs})d\theta.
\end{eqnarray*}
Furthermore, denote the stationary distributions of the Gibbs chain and the iterative chain by $\nu_1^{\xx^{obs}}$ and $\nu_2^{\xx^{obs}}$.

\begin{theorem}\label{PropNece}
Suppose that for some $j\in \mathbb Z^+$, sets $A$ and $C$, and $\varepsilon\in(0,1/2)$
\begin{equation*}
\inf_{\xx^{mis}_{-j}\in C}P^{g}_j(\xx_j^{mis}\in A|\xx^{mis}_{-j},\xx^{obs})>
\sup_{\xx^{mis}_{-j}\in C}P^{f}(\xx_j^{mis}\in A|\xx^{mis}_{-j},\xx^{obs})+\varepsilon
\end{equation*} or
\begin{equation*}
\sup_{\xx^{mis}_{-j}\in C}P^{g}_j(\xx_j^{mis}\in A|\xx^{mis}_{-j},\xx^{obs})<
\inf_{\xx^{mis}_{-j}\in C}P^{f}(\xx_j^{mis}\in A|\xx^{mis}_{-j},\xx^{obs})-\varepsilon
\end{equation*}
and $\nu_1^{\xx^{obs}}(\xx_{-j}^{mis}\in
C)>q \in (0,1)$.

Then there exists a set $B$ such that
\begin{equation*}
\left|\nu_2^{\xx^{obs}}(\xx^{mis}\in B) - \nu_1^{\xx^{obs}}(\xx^{mis}\in
B)\right|> q \varepsilon /4.
\end{equation*}

\begin{proof}
Suppose that \begin{equation*} \inf_{\xx^{mis}_{-j}\in C}P^{g}_j(\xx_j^{mis}\in A|\xx^{mis}_{-j},\xx^{obs})>
\sup_{\xx^{mis}_{-j}\in C}P^{f}(\xx_j^{mis}\in A|\xx^{mis}_{-j},\xx^{obs})+\varepsilon,
\end{equation*}
The ``less than'' case is completely analogous. Consider
the set $B=\{\xx^{mis}: \xx_{-j}^{mis}\in C,
\xx_{j}^{mis}\in A\}$. If
\begin{equation}\label{CC}|\nu_2^{\xx^{obs}}(\xx^{mis}_{-j}\in C) -
\nu_1^{\xx^{obs}}(\xx^{mis}_{-j}\in C)| \leq q\varepsilon/2 ,
\end{equation}then, by the fact that
\begin{eqnarray*}\nu_1^{\xx^{obs}}(\xx^{mis}\in B)=\nu_1^{\xx^{obs}}(\xx^{mis}_{-j}\in C)\int P^{f}(\xx_j^{mis}\in
A|\xx^{mis}_{-j},\xx^{obs}) \nu_1^{\xx^{obs}}(d\xx^{mis}_{-j}|\xx^{mis}_{-j}\in C),\\
\nu_2^{\xx^{obs}}(\xx^{mis}\in B)=\nu_2^{\xx^{obs}}(\xx^{mis}_{-j}\in C)\int P^{g}(\xx_j^{mis}\in
A|\xx^{mis}_{-j},\xx^{obs}) \nu_2^{\xx^{obs}}(d\xx^{mis}_{-j}|\xx^{mis}_{-j}\in C),
\end{eqnarray*}
we obtain
\begin{equation*}
|\nu_2^{\xx^{obs}}(\xx^{mis}\in B) - \nu_1^{\xx^{obs}}(\xx^{mis}\in
B)|> q \varepsilon /4.
\end{equation*} Otherwise, if \eqref{CC} does not hold, let $B=\{\xx^{mis}: \xx^{mis}_{-j}\in
C\}$.
\end{proof}
\end{theorem}
For two models with different likelihood functions, one can construct sets $A$ and $C$ such that the conditions in the above theorem hold. Therefore, if among the predictive distributions of all the $p$ conditional
models there is one $g_j$ that is different from $f$ as stated in Theorem \ref{PropNece}, then the stationary distribution of the
iterative imputation is different from the
posterior distribution of the Bayesian model in total variation
by a fixed amount. For a set of incompatible models and any
joint model $f$, there exists at least one $j$ such that the
conditional likelihood functions of $\xx_j$ given $\xx_{-j}$
are different for $f$ and $g_j$. Their predictive distributions have to be
different for $\xx_j$. Therefore, such an iterative
imputation using incompatible conditional models typically does not correspond to Bayesian imputation under any joint model.

\section{Incompatible conditional models}

\label{SecIncomp}

In this section, we proceed to the discussion of incompatible conditional models. We first extend the concept of model compatibility to semi-compatibility which includes the regression models we have generally seen in practical uses of iterative imputation. We then introduce the validity of semi-compatible models. Finally, we show that if the conditional models are semi-compatible and valid (together with a few mild technical conditions) the combined imputation estimator is consistent.

\subsection{Semi-compatibility and model validity}

As in the previous section, we assume that the invariant distribution  exists. For compatible conditional models, we used the posterior distribution of the corresponding Bayesian model as the natural benchmark and show that the two
imputation distributions converge to each other. We can use this idea for the analysis of incompatible models. In this setting, the first issue
is to find a natural Bayesian model associated with a set of incompatible
conditional models. Naturally, we introduce the concept of semi-compatibility.

\begin{definition}
\label{defSemiComp} A set of conditional models $\{h_j(x_j|x_{-j},\theta_j,%
\varphi_j): j=1,...,p\}$, each of which is indexed by two sets of parameters
$(\theta_j,\varphi_j)$, is said to be \emph{semi-compatible}, if there
exists a set of compatible conditional models
\begin{equation}  \label{SemiComp}
g_j(x_j|x_{-j},\theta_j)=h_j(x_j|x_{-j},\theta_j,\varphi_j=0),
\end{equation}
for $j=1,...,p$. We call $\{g_j:j=1,...,p\}$ a \emph{compatible element} of $%
\{h_j: j=1,...,p\}$.
\end{definition}

By definition, every set of compatible conditional models is
semi-compatible.  A simple and uninteresting class of semi-compatible models arises with iterative regression imputation.  As typically parameterized, these models include complete independence as a special case. A \emph{trivial} compatible element, then, is the one in which $x_j$ is independent of $x_{-j}$ under $g_j$ for all $j$.  Throughout the discussion of this section, we use $\{g_j: j=1,...,p\}$ to denote the compatible element of $\{h_j:j=1,...,p\}$
and $f$ to denote the joint model compatible with $\{g_j:j=1,...,p\}$.

Semi-compatibility is a natural concept connecting a joint
probability model to a set of conditionals.
One foundation of almost all statistical theories is that data are
generated according to some (unknown) probability law.
When setting up each conditional model, the
imputer chooses a rich family that is intended to include distributions that are close to the true conditional distribution. For instance, as recommended by \cite{Meng1994}, the imputer should try to
include as many predictors as possible (using regularization as necessary to keep the estimates stable). Sometimes, the degrees of
flexibility among the conditional models are different. For instance, some
includes quadratic terms or interactions. This richness usually results in
incompatibility. Semi-compatibility includes such cases in which the
conditional models are rich enough to include the true model but may not be
always compatible among themselves. To proceed, we introduce the following
definition.

\begin{definition}
\label{DefValid} Let $\{h_j:j=1,...,p\}$ be semi-compatible, $\{g_j:j=1,...,p\}$ be its compatible element, and $f$ be the joint model compatible with $g_j$. If the joint model $f(x|\theta)$ includes the true probability distribution, we say $\{h_j:j=1,...,p\}$ is a set of \emph{valid
semi-compatible models}.
\end{definition}

In order to obtain good prediction, we need the validity of the semi-compatible models. A natural
issue is the performance of valid semi-compatible models. Given that we have
given up compatibility, we should not expect the iterative imputation to be equivalent to any joint Bayesian imputation. Nevertheless, under mild conditions, we
are able to show the consistency of the combined imputation estimator.

\subsection{Main theorem of incompatible conditional models}

 Now, we list a set of  conditions:

\begin{enumerate}
\item [B1] Both the Gibbs and iterative chains admit their unique invariant distributions, $\nu_1^{\xx^{obs}}$ and $\nu_2^{\xx^{obs}}$.

\item[B2] The posterior distributions of $\theta $ (based on $f$) and $(\theta_j, \varphi_j)$ (based on $h_j$) given a complete data set $\xx$ have the representation
$|\theta -\tilde  \theta| \leq \xi n^{-1/2}, |(\theta_{j}-\tilde \theta_{j},\varphi_j - \tilde \varphi_j)|\leq \xi_{j} n^{-1/2},$
where $\tilde \theta$ is the maximum likelihood estimate of $f(\xx|\theta)$, $(\tilde \theta_j, \tilde \varphi_j)$ is the maximum likelihood estimate of $h_j(\xx_j|\xx_{j},\theta_j,\varphi_j)$, and $Ee^{|\xi_{j}|}\leq \kappa $, $E e^{|\xi|} \leq \kappa$ for some $\kappa >0$.

\item [B3] All the score functions  have finite moment generating functions under $f(\xx^{mis} |\xx^{obs}, \theta)$.

\item[B4] For each variable $j$, there exists a subset of observations $\iota_j$ so that for each  $i \in\iota_{j}$  $x_{i,j}$ is missing and $x_{i,-j}$ is fully observed. In addition, the cardinality $\#(\iota_j) \rightarrow \infty$ as $n\rightarrow \infty$.
\end{enumerate}

\begin{remark}
Conditions B2 and B3 impose moment conditions on the posterior distribution and the score functions. They are satisfied by most parametric families. Condition B4 rules out certain boundary cases of missingness patterns and is imposed for technical purposes. The condition it is not very restrictive because it only requires that the cardinality of $\iota_{j}$ tends to infinity, not necessarily even of order of $O(n)$.
\end{remark}

We now express the fifth and final condition which requires the following construction. 
Assume the conditional models are valid and that the data $\xx$ is generated from $f(\xx|\theta^0)$. We use $\theta_{j}^{0}=t_{j}(\theta ^{0})$ and $\varphi _{j}^{0}=0$ to
denote the true parameters under $h_{j}$. We define
\begin{equation}\label{EE}
\hat{\theta}=\sup_{\theta }f(\xx^{obs}|\theta ),\qquad\hat{\theta}_{j}=t_{j}(\hat{%
\theta}),
\end{equation}
be the observed-data MLE and
\begin{equation}\label{note}
\hat\theta^{(2)} = \arg\sup_\theta \int \log  f(\xx|\theta) \nu^{\xx^{obs}}_2 (d\xx^{mis}),\quad
(\hat\theta_j^{(2)},\hat\varphi_j^{(2)}) = \arg\sup_{\theta_j,\varphi_j} \int  \log h_j(\xx_j|\xx_{-j},\theta_j,\varphi_j) \nu^{\xx^{obs}}_2 (d\xx^{mis})
\end{equation}
where $\xx= (\xx^{obs},\xx^{mis})$.

Consider a Markov chain $x^{\ast }(k)$ corresponding to one observation---one row of the data matrix---living on $R^{p}$. The chain evolves as follows. Within each iteration, each dimension $j$ is
updated conditional on the others according to the conditional distribution
\begin{equation*}
h_{j}(x_{j}|x_{-j},\theta _{j},\varphi _{j}),
\end{equation*}%
where $(\theta _{j},\varphi _{j})=(\hat{\theta}_{j},0)+\varepsilon \xi _{j}$ and $\xi_j$ is a random vector with finite MGF (independent of everything at every step).
Alternatively, one may consider $(\theta_{j},\varphi_{j})$ as a sample from the posterior distribution corresponding to the conditional model $h_{j}$.
 Thus, $x^*(k)$ is the marginal chain of one observation in iterative chain. Given that $\xx^{mis,2}(k)$ admits a unique invariant distribution, $x^{\ast }(k)$ admits its unique stationary distribution  for $\varepsilon$ sufficiently small.
Furthermore, consider that $x(k)$  is a Gibbs sampler and it admits stationary distribution  $f(x|\hat{\theta})$, that is, each component is updated according to the conditional distribution $f(x_j|x_{-j},\hat \theta)$ and the parameters of the updating distribution are set at the observed data maximum likelihood estimate, $\hat \theta$. The last condition is stated as follows.

\begin{enumerate}
\item[B5]
$x^{\ast }(k)$ and $x(k)$ satisfy conditions in Lemma \ref
{ThmComp} as $\varepsilon \rightarrow 0$, that is, the invariant distributions of $x^{*}(k)$ and $x(k)$ converges in $\Vert\cdot \Vert_{V}$ norm, where $V$ is a drift function for $x^*(k)$. There exists a constant $\kappa$ such that all the score functions are bounded by 
$$\partial\log f(x| \theta^{0})\leq \kappa V(x),\qquad  \partial\log h_j(x_j|x_{-j},t_{j}(\theta^{0}),\varphi_j=0)\leq \kappa V(x) .$$
\end{enumerate}

\begin{remark}
By choosing $\varepsilon$ small, the transition kernels of $x^*(k)$ and $x(k)$ converge to each other. Condition B5 requires that Lemma \ref{ThmComp} applies in this setting, that their invariant distributions are close in the sense stated in the Lemma. This condition does not suggest that Lemma \ref{ThmComp} applies to $\nu^{\xx^{obs}}_1$ and $\nu^{\xx^{obs}}_2$, which represents the joint distribution of many such $x^*(k)$'s and $x(k)$'s.
\end{remark}

We can now state the main theorem in this section.

%

\begin{theorem}\label{ThmIncomp} Consider a set of valid semi-compatible models $\{h_{j}:j=1,...,p\}$, and assume conditions B1--5 are in force.  Then, following the notations in \eqref{note}, the following limits hold:
\begin{equation}\label{consist}
\hat\theta^{(2)}\rightarrow \theta ^{0},\qquad
\hat\theta_j^{(2)}\rightarrow t_j(\theta^0),
\quad \hat\varphi_j^{(2)}\rightarrow 0,
\end{equation}
in probability as sample size $n\rightarrow \infty$ for all $j$.
\end{theorem}

\begin{remark}
The expression $\hat \theta^{(2)}$ corresponds to the following estimator. 
Impute the missing data from distribution $\nu^{\xx^{obs}}_{2}$ $m$ times to obtain $m$ complete datasets. Stack the $m$ datasets to one big dataset. Let $\hat \theta^{(2)}_{m}$ be the maximum likelihood estimator based on the big dataset. Then, $\hat \theta^{(2)}_{m}$ converges to $\hat \theta^{(2)}$ as $m\rightarrow \infty$.
Furthermore, $\hat \theta^{(2)}$ is asymptotically equivalent to 
the combined point estimator of $\theta$ according to Rubin's combining rule (with infinitely many imputations). Similarly, $(\hat \theta_j^{(2)}, \hat \varphi_j^{(2)})$ is asymptotically equivalent to the combined estimator of the conditional model. Therefore, Theorem \ref{ThmIncomp} suggests that the combined imputation estimators are consistent under conditions B1--5.
\end{remark}


\section{Linear example}\label{SecLinear}

\subsection{A simple set of compatible conditional models}

\begin{figure}
\begin{center}
\includegraphics[height=2in]{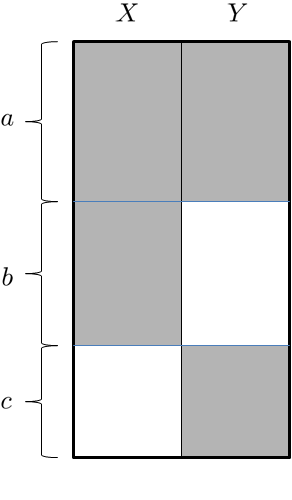}
\end{center}
\caption{Missingness pattern for our simple example with two variables.  Gray and white areas indicate observed and missing data, respectively.   This example is constructed so that there are no cases for which both variables are missing.}\label{FigData}
\end{figure}

In this subsection, we study a linear model as an illustration. Consider $n$ i.i.d.~bivariate
observations $(\xx,\yy)=\{(x_i,y_i): i=1,...,n\}$ and a set of
conditional models
\begin{equation}\label{IterLinear}x_i | y_i
\sim N(\beta_{x|y} y_i, \tau_x^2), \quad y_i |x_i \sim N(\beta_{y|x}
x_i, \tau_y^2).\end{equation}
To simplify the discussion, we
set the intercepts to zero. As discussed previously, the joint
compatible model assumes that $(x,y)$ is a bivariate normal
random variable with mean zero, variances $\sigma_x^2$ and
$\sigma_y^2$, and correlation $\rho$. The reparameterization
from the joint model to the conditional model of $y$ given $x$
is
$$\beta_{y|x}= \frac{\sigma_y}{\sigma_x}\rho, \quad \tau_y^2 = (1-\rho^2) \sigma_y^2.
$$
Figure \ref{FigData} displays the missingness pattern we are assuming for this simple example, with $a$ denoting the set of
observations for which both $x$ and $y$ are observed, $b$ denote those with
missing $y$'s, and $c$ denoting those with missing $x$'s; $n_a$, $n_b$,
and $n_c$ denote their respective sample sizes, and $n= n_a + n_b + n_c$. To keep the example simple, we assume that there are no cases for which both $x$ and $y$ are missing.

\paragraph{Positive recurrence and limiting distributions.}

The Gibbs chain and the iterative chain admit a common small set containing the observed-data maximum likelihood estimate. The construction of the drift functions is tedious and is not particularly relevant to the current discussion, and so we leave their detailed derivations to the supplemental materials available at http://stat.columbia.edu/$\sim$jcliu/paper/driftsupp.pdf. We proceed here by assuming that they are in force.

\paragraph{Total variation distance between the kernels.}
The results for incompatible models apply here. Thus, condition A1 in Theorem \ref{ThmCompRev} has been satisfied.
We now check the boundedness of $\partial L(\theta)$. The posterior distribution of the full Bayes
model is
\begin{eqnarray*}
p(\sigma_x^2,\tau_y^2,\beta_{y|x}|\xx,\yy)&\propto& f( \xx,\yy|\sigma_x^2,\tau_y^2, \beta_{y|x})
\pi^*(\sigma_x^2,\tau_y^2, \beta_{y|x}) \\
&=&f(\yy|\tau_y^2,\beta_{y|x},\xx)
f(\xx | \sigma_x^2)\pi^*(\sigma_x^2,\tau_y^2, \beta_{y|x}).
\end{eqnarray*}
The posterior distribution of $(\tau_y^2,\beta_{y|x})$ with $\sigma_x^2$ integrated out is
\begin{equation*}
p(\tau_y^2,\beta_{y|x}|\xx,\yy)\propto f(\yy%
|\tau_y^2,\beta_{y|x},\xx)\pi_{\xx}(\beta_{y|x},\tau_y^2),
\end{equation*}
where
\begin{equation*}
\pi_{\xx}(\beta_{y|x},\tau_y^2) \propto \int f(\xx|
\sigma_x^2)\pi^*(\sigma_x^2,\tau_y^2, \beta_{y|x}) d\sigma_x^2.
\end{equation*}
The next task is to show that $\pi_{\xx}(\beta_{y|x},\tau_y^2)$ is a diffuse prior satisfying the conditions in Proposition \ref{PropPrior}. We impose the following independent prior distributions on $\sigma_x^2$, $\sigma_y^2$, and $\rho$:
\begin{equation}\label{prior}
\pi(\sigma_x^2,\sigma_y^2,\rho )\propto \sigma_x\sigma_y I_{[-1,1]}(\rho).
\end{equation}
The distribution of $\xx$ does not depend on $(\sigma_y^2,\rho)$. Therefore, under the posterior distribution given $\xx$, $\sigma_x^2$ and $(\sigma_y^2,\rho)$ are independent.  Conditional on $\xx $, $\sigma_x^2$ is inverse-gamma. Now we proceed to develop the conditional/posterior distribution of $(\tau_y^2, \beta_{y|x})$ given $\xx$. Consider the following change of variables
\begin{equation*}
\sigma_y^2 = \tau_y^2 + \beta_{y|x}^2 \sigma_x^2, \quad \rho = \beta_{y|x}\sqrt{\frac{
\sigma_x^2}{\tau_y^2 + \beta_{y|x}^2 \sigma_x^2}}.
\end{equation*}%
Then,
\begin{equation*}
\det \left(\frac{\partial (\sigma_y^2,\rho,\sigma_x^2)}{\partial
(\tau_y^2, \beta_{y|x},\sigma_x^2)}\right)= \frac{\sigma_x }{\sqrt{\tau_y^2 + \beta_{y|x}^2 \sigma_x^2}}.
\end{equation*}
Together with
\begin{equation*}
\pi(\sigma_y^2, \rho^2)\propto \sigma_y^{},
\end{equation*}
we have
\begin{eqnarray*}
\pi_{\xx}(\tau_y^2, \beta_{y|x})&\propto& \int \det \left(\frac{\partial (\sigma_y^2,\rho,\sigma_x^2)}{\partial
(\tau_y^2, \beta_{y|x},\sigma_x^2)}\right) \pi(\sigma^2_y, \rho) p(\sigma_x^2|\xx )d\sigma_x^2 \\&=&\int \sigma_x p(\sigma_x^2|\xx )d\sigma_x^2= C(\xx).
\end{eqnarray*}
\begin{remark}
If one chooses $\pi_2 (\tau_y^2,\beta_{y|x})\propto 1$ for the iterative imputation and (\ref{prior}) for the joint Bayesian model, the iterative chain and the Gibbs chain happen to have identical transition kernels and, therefore, identical invariant distributions. This is one of the rare occasions that these two procedures yield identical imputation distributions.
\end{remark}

If one chooses Jeffreys' prior, $\pi_2 (\tau_y^2,\beta_{y|x})\propto \tau_y^{-2}$, then
$$L(\tau_y^2,\beta_{y|x})=\frac{\pi_{\xx}(\tau^2_y,\beta_{y|x})}{\pi_2 (\tau_y^2,\beta_{y|x})}\propto \tau_y^2,$$
and $\partial L$ is bounded in a suitably chosen compact set containing the true parameters. Thus, Theorem \ref{ThmCompRev} applies.


\paragraph{Empirical convergence check.}


\begin{figure}
\begin{center}
\includegraphics[height=2in]{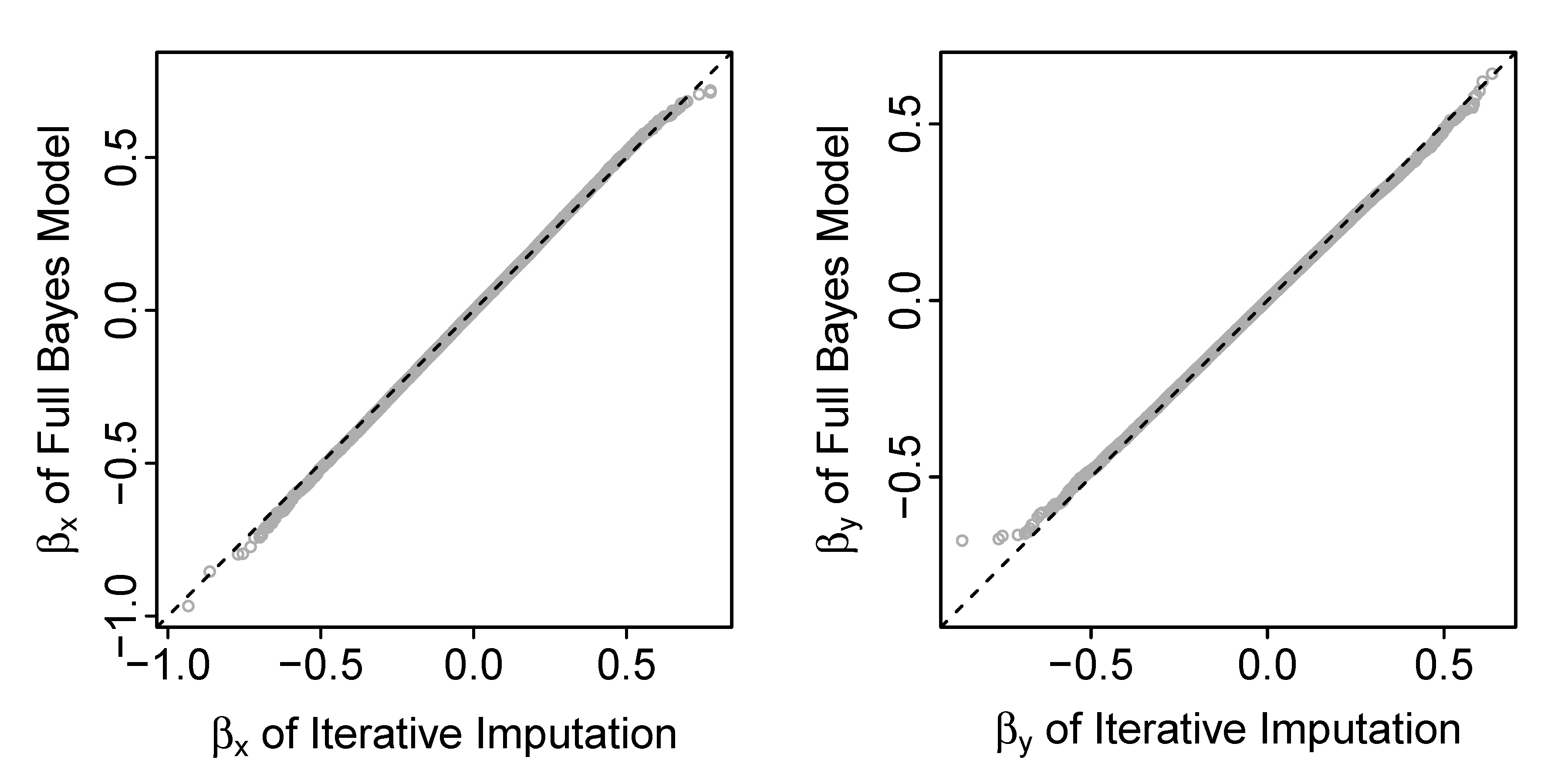}
\end{center}
\caption{Quantile-quantile plots demonstrating the closeness of the posterior distribution of the Bayesian model and the compatible iterative
imputation distributions for $\beta_{x}$ and $\beta_{y}$ with sample size $n_a=200$.}\label{FigPost200}
\end{figure}



To numerically confirm the convergence of the two distributions, we
generate the following data sets. To simplify analysis, let $(x_i,y_i)$'s be bivariate Gaussian random vectors with mean zero, variance one, and correlation zero. We set $n_a = 200$, $n_b = 80$, and $n_c = 80$. For the iterative imputation we use Jeffreys' prior $p(\tau_y^2,\beta_{y|x})\propto \tau_y^{-2}$ and $p(\tau_x^2,\beta_{x|y})\propto \tau_x^{-2}$. For the full Bayesian model, the prior distribution is chosen as in (\ref{prior}).

We monitor the posterior distributions of the following statistics:
\begin{equation}\label{monitor}\beta_{x}= \frac{\sum_{i\in b}x_iy_i}{\sum_{i\in
b}y_i^2},\quad \beta_{y}= \frac{\sum_{i\in
c}x_iy_i}{\sum_{i\in c}x_i^2}.\end{equation}
Figures \ref{FigPost200} shows the quantile-quantile plots of the distributions of $\beta_x$ and $\beta_y$ under $\nu_1^{\xx^{obs}}$ and $\nu_2^{\xx^{obs}}$ based on $1$ million MCMC iterations. The differences between these two distributions are tiny.

\subsection{Higher-dimensional linear models}

We next consider a more complicated and realistic situation, in which there are $p$ continuous variables, $x_{1},...,x_{p}$. Each conditional model is linear in the sense that, for each $j$,
\begin{equation*}
x_{j}|x_{-j}\sim N((1, x_{-j}^{\top})\beta_{j},\sigma_{j}^{2}),
\end{equation*}
which is the set of compatible models presented in Example \ref{ExLinear}.

\begin{figure}
\centering
  \includegraphics[width=\textwidth]{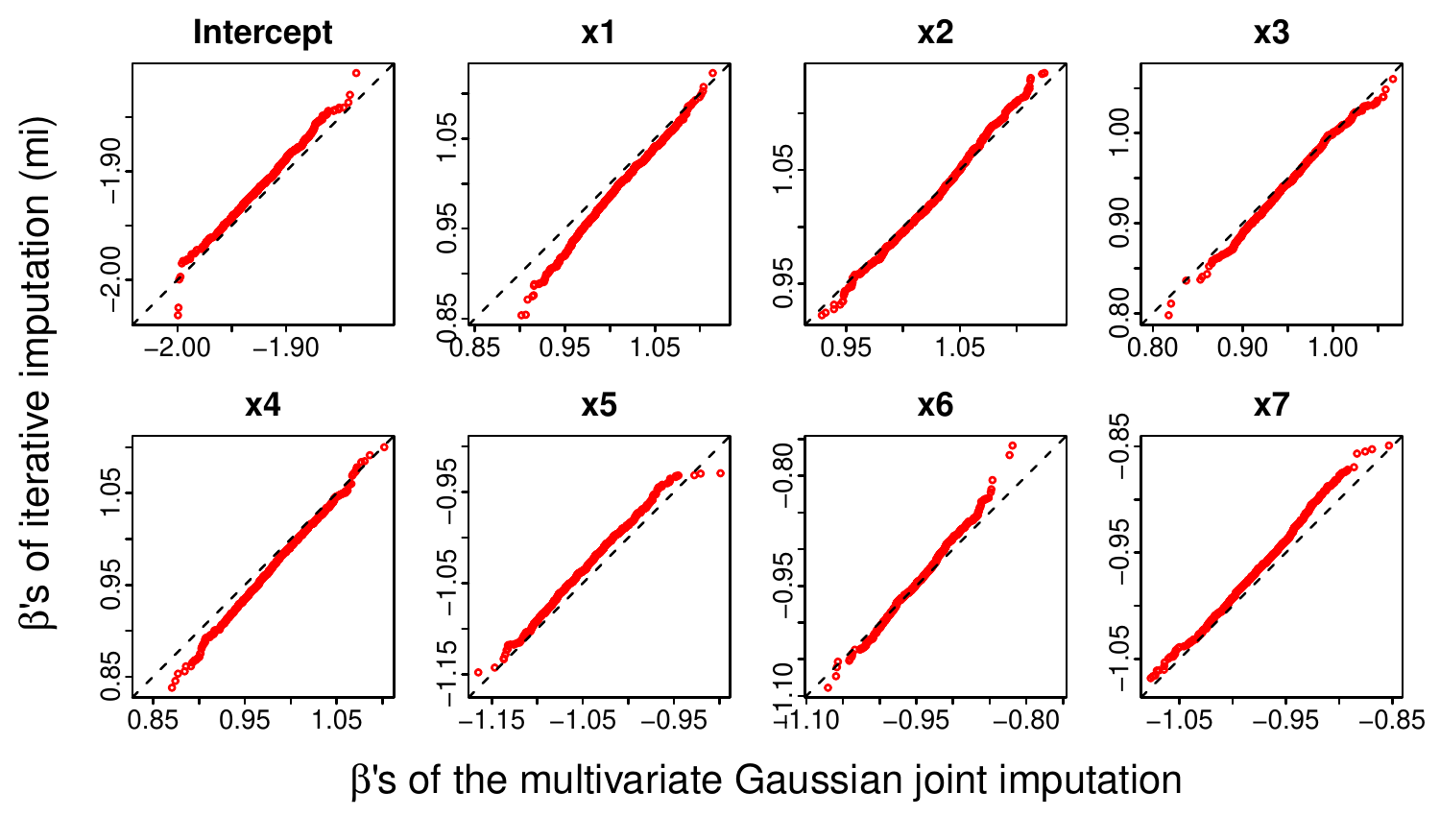}\\
  \caption{Quantile-quantile plots of the imputation distributions of the regression coefficients ($y$ on $x$'s) from  the joint Bayesian imputation and the iterative imputation.}\label{fig:normbymi}
\end{figure}

In the simulation, we generate 1000 samples of $(x_{1},...,x_7)$ from a 7-dimensional multivariate normal distribution with mean 0 and covariance matrix that equals 1 on the diagonals and 0.4 on the off-diagonal elements.
We then generate another variable $y\sim N( -2 + x_1 + x_2 + x_3 + x_4 - x_5 - x_6 - x_7, 1)$. 
Hence the  dataset contains $y, x_1, x_2, \ldots, x_7$.
For each variable, we randomly select $30\%$ of the observations and set them to be missing. Thus, the missing pattern of the dataset is missing completely at random (MCAR). We impute the missing values in two ways:  iterative imputation and a multivariate Gaussian joint Bayesian model.
After imputation, we use the imputed datasets and regress $y$ on all $x$'s to obtain the regression coefficients. The quantile-quantile plots in Figure \ref{fig:normbymi} compare the imputation distribution of the least-square estimates of the regression coefficients  of  the iterative imputation procedure and the multivariate Gaussian joint model. 


\subsection{Simulation study for incompatible models}

We next consider conditional models that are incompatible and valid.
To study the frequency properties of the iterative imputation algorithm, we generate 1000  datasets independently each with a sample size of 10,000.  For each dataset, $y_1\sim\mbox{Bernouli}(0.45)$, $y_2\sim\mbox{Bernouli}(0.65)$, $y_{1}$ and $y_{2}$ are independent, and the remaining variables come from this conditional distribution: $x_1,\ldots, x_5|y_1,y_2\sim N(\mu_1 y _1 + \mu_2 y_2, \Sigma)$, where $\mu_1$ is a vector of 1's and $\mu_2$ is a vector of $0.5$'s and $\Sigma$ is a $5\times 5$ matrix that is 1 on the diagonals and 0.2 on the off-diagonal elements.

\begin{figure}
\centering
  \includegraphics[width=\textwidth]{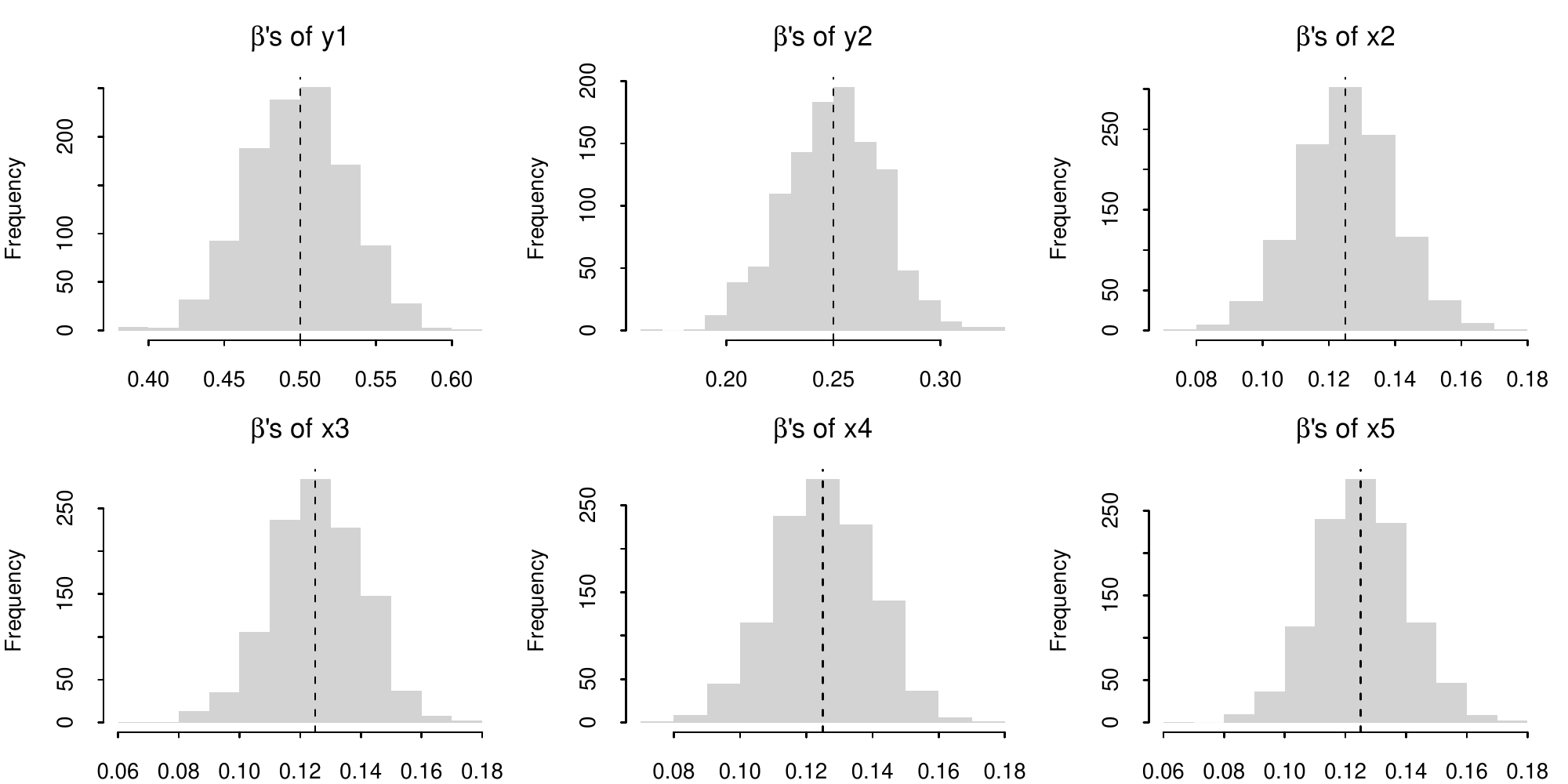}\\
  \caption{Distributions of coefficients of $x_1$ regressing on $y_1, y_2, x_2, x_3, x_4, x_5$ from 1000 imputed datasets using an iterative imputation routine \cite{Su2011}. The dashed vertical lines represents the true value of the regression coefficients of the simulation setting, which are $0.5, 0.25, 0.125, 0.125, 0.125, 0.125$.}\label{fig:histglm_hi}
\end{figure}

As before, we remove 30\% of the data completely at random and then impute the dataset using iterative imputation.  We impute $y_1$ and $y_2$ using logistic regressions and $x_1, \ldots, x_5$ using linear regressions. In particular, $y_1$ is conditionally imputed given $y_{2}, x_{1}, x_{2}, x_{3}, x_{4}, x_{5}$, and the interactions $x_1y_2$ and $x_2y_2$;
$y_2$ is conditionally imputed given $y_{1}, x_{1}, x_{2}, x_{3}, x_{4}, x_{5}$, and the interactions $x_1y_1$ and $x_2y_1$; and each $x_j$, $j=1,\dots,5$, is conditionally imputed given $y_1,y_2,$ and the other four $x_j$'s.
The conditional models for the $x_{j}$'s are simple linear models, whereas the logistic regressions for $y_{i}$ also include interactions.  As a result, the set of conditional models is no longer compatible but is still valid.
To check whether or not the incompatible models result in reasonable estimates,  we impute the missing values using these conditional models. For each dataset, we obtain combined estimates of the regression coefficients of $x_{1}$ given the others by averaging the least-square estimates over 50 imputed datasets.
That is, for each dataset, we have 50 imputations, for each of which we obtain the estimated  regression coefficients of $x_1|y_1, y_2, x_2, x_3, x_4, x_5$.  Next, we average over 50 sets of coefficients to obtain a single set of coefficients.  We repeat the whole procedure on 1000 datasets to get 1000 sets of estimated coefficients. 
Figure \ref{fig:histglm_hi} shows the distribution of the estimated coefficients of $x_1$ regressing on $y_1, y_2, x_2, x_3, x_4, x_5$ based on  the 1000 independent datasets.  The frequentist distributions of the combined estimate are centered around their true values indicated by the dashed line.
This is consistent with Theorem \ref{ThmIncomp}.

\section{Discussion}

Iterative imputation is appealing in that it promises to solve the difficult task of multivariate modeling and imputation using the flexible and simple tools of regression modeling.  But two key concerns arise:  does the algorithm converge to a stationary distribution and, if it converges, how to interpret the resulting joint distribution of the imputations, given that in general it will not correspond exactly to the fitted conditionals from the regression.
In this article, we have taken steps in that direction.


There are several natural directions for future research.  From one direction, it should be possible to obtain exact results for some particular classes of models such as linear regressions with Gaussian errors and Gaussian prior distributions, in which case convergence can be expressed in terms of simple matrix operations.  In the more general case of arbitrary families of regression models, it would be desirable to develop diagnostics for stationarity (along with proofs of the effectiveness of such diagnostics, under some conditions) and empirical measures of the magnitude of discrepancies between fitted and stationary conditional distributions.

Another open problem here is how to consistently estimate the variance of the combined imputation estimator. Given that the imputation distribution of incompatible models is asymptotically different from that of any joint Bayesian imputation, there is no guarantee that Rubin's combined variance estimator is asymptotically consistent. We acknowledge that this is a challenging problem. Even for joint Bayesian imputation, estimating the variance of the combined estimator is still a nontrivial task under specific situations; see, for instance, \cite{Kim04,Meng1994}. Therefore, we leave this issue to future studies.

We conclude with some brief notes.

\paragraph{A special case of the compatible models.}
In the analysis of the conditional models, suppose that the parameter spaces of the conditional distribution and the covariates are separable, that is, $f(x_{j},x_{-j}|\theta_{j},\theta_{j}^{*}) = f(x_{j}|x_{-j},\theta_{j}) f(x_{-j}|\theta_{j}^{*})$ and there exists a prior $\pi$ for the joint model $f$ such that  $\theta_{j}$ and $\theta_{j}^{*}$ are a priori independent for all $j$.
The,  the boundedness of $\partial\log L({\theta_{j}})$ becomes straightforward to obtain.  Note that $L(\theta_{j}) = \pi(\theta_{j})/ \pi_{j,\xx_{-j}}(\theta_{j})$ and
\begin{equation*}
\pi_{j,\xx_{-j}}(\theta_{j})\triangleq \pi_{j}^{\ast }(\theta_{j})  \int f(%
\xx_{-j}|,\theta_{j}^{\ast })\pi_{j}^{\ast }(\theta_{j}^{*})d\theta_{j}^{\ast }.
\end{equation*}%
Thus, $L(\theta_{j}) = \pi(\theta_{j})/ \pi_{j}^{*}(\theta_{j})$ is independent of the data.

\paragraph{A example when Theorem \ref{ThmCompRev} does not apply.}
The expression $\partial\log L(\theta_{j})$ in Proposition \ref{PropPrior} is not always bounded. For example, suppose that $\pi_{j,\xx_{-j}}(\theta_{j})$ is an informative prior for which the covariates of the regression model $\xx_{-j}$ provide strong information on the regression coefficients $\theta_{j}$ according to the joint model $f$. For instance,  in  Example \ref{ExLogit} on page \pageref{ExLogit}, the marginal distribution of the covariate $x_{2}$ provides $\sqrt n$ amount of information on the coefficients of the logistic regression. Under this situation, Proposition \ref{PropPrior} still holds, but the right-hand side of \eqref{pbd} does not necessarily converge to zero. Consequently, Theorem \ref{ThmCompRev} (the main result in this section, presented later) does not apply. For the properties of iterative imputation under such situations, we can apply the consistency result for incompatible models, as discussed  in Section \ref{SecIncomp}. That is, the combined estimator is still consistent though it is not equivalent to any Bayesian model.



\appendix

\section{Proofs in Section \ref{SecCompatible}} \label{Apdbd}

\begin{lemma}
\label{LemTV1}Let $Q_{0}$ and $Q_{1}$ be probability measures defined on the
same $\sigma$-field $\mathcal{F}$ and such that $dQ_{1}=r^{-1}dQ_{0}$ for a
positive r.v. $r>0$. Suppose that for some $\varepsilon>0$, $E^{Q_{1}}\left(
r^{2}\right)  =E^{Q_{0}}r\leq1+\varepsilon$. Then,
\[
\sup_{|f|\leq 1}\left\vert E ^{Q_{1}}(  f(X))  - E^{Q_{0}}(
f(X))  \right\vert \leq\varepsilon^{1/2}.
\]
\end{lemma}

\begin{proof}[Proof of Lemma \ref{LemTV1}]
\begin{align*}
\left\vert E ^{Q_{1}}(  f(X))  - E^{Q_{0}}(
f(X))  \right\vert  &
=\left\vert E^{Q_{1}}\left[  (1-r) f(X)\right]  \right\vert \\
&  \leq E^{Q_{1}}\left(  \left\vert r-1\right\vert \right)  \leq [E^{Q_{1}}(  r-1)  ^{2}]^{1/2}=\left(  E^{Q_{1}}r^{2}-1\right)  ^{1/2}%
\leq\varepsilon^{1/2}.
\end{align*}
\end{proof}

\bigskip

\begin{proof}[Proof of Proposition \ref{PropPrior}]
From Lemma \ref{LemTV1}, we need to show that
\begin{equation*}
\int r^{2}(\theta )f_{\xx}(\theta )d\theta \leq 1+\frac{\kappa^2|\partial L(\mu _{\theta })|}{\sqrt n L(\mu
_{\theta })}.
\end{equation*}%
Let $\mu _{L}=E^{f}L(\theta )$.
\begin{equation*}
r(\theta )=\frac{L(\theta )}{\mu _{L}}=\frac{L(\mu _{\theta })+\partial
L(\xi )(\theta -\mu _{\theta })}{\mu _{L}}.
\end{equation*}%
Then%
\begin{eqnarray*}
E^{f}(r^{2}(\theta ))\frac{\mu _{L}^{2}}{L^{2}(\mu _{\theta })} &=&1+2E\frac{%
\partial L(\xi )(\theta -\mu _{\theta })}{L(\mu _{\theta })}+E\frac{(\partial
L(\xi ))^{2}(\theta -\mu _{\theta })^{2}}{L^{2}(\mu _{\theta })} \\
&\leq &1+2\frac{|\partial L(\mu_{\theta} )|E Z}{L(\mu _{\theta })\sqrt{n}}+\frac{%
(\partial L(\mu _{\theta }))^{2}}{L^{2}(\mu _{\theta })}\frac{EZ^{2}}{n}.
\end{eqnarray*}%
With a similar argument, there exists a constant $\kappa_{1}$ such that%
\begin{equation*}
\left\vert \frac{\mu _{L}^{2}}{L^{2}(\mu _{\theta })}-1\right\vert \leq
\frac{|\partial L(\mu _{\theta })|}{L(\mu _{\theta })}\frac{\kappa^2_1 }{\sqrt{n}}.
\end{equation*}%
Therefore, there exists some $\kappa>0$ such that%
\begin{eqnarray*}
E^{f}r^{2}(\theta ) &\leq &\left( 1+2E\frac{|\partial L(\mu_{\theta} )|Z}{L(\mu
_{\theta })\sqrt{n}}+\frac{(\partial L(\mu _{\theta }))^{2}}{L^{2}(\mu
_{\theta })}\frac{EZ^{2}}{n}\right) \frac{L^{2}(\mu _{\theta })}{\mu _{L}^{2}%
} \\
&\leq &1+\frac{|\partial L(\mu _{\theta })|}{L(\mu _{\theta })}\frac{%
\kappa }{\sqrt{n}}.
\end{eqnarray*}%
Using Lemma \ref{LemTV1}, we conclude the proof.
\end{proof}

\bigskip

\begin{proof}[Proof of Lemma \protect\ref{ThmComp}]
For any $\varepsilon ,\delta >0$, let $k_{\varepsilon }=\inf \{j:\forall
k>j,r_{k}\leq \varepsilon \}$. Then, for any $m>k_{\varepsilon }$
\begin{eqnarray}
\Vert\tilde \nu^{\xx^{obs}} _{1}-\tilde \nu^{\xx^{obs}} _{2}\Vert_{V} &\leq &\left\Vert \tilde \nu^{\xx^{obs}}
_{1}-\frac{1}{m}\sum_{k=1}^{m}\tilde K_{1}^{(k)}(\nu ,\cdot
)\right\Vert _{V}+\left\Vert \tilde \nu^{\xx^{obs}} _{2}-\frac{1}{m}
\sum_{k=1}^{m}\tilde K_{2}^{(k)}(\nu ,\cdot )\right\Vert _{V}  \notag \\
&&+\left\Vert \frac{1}{m}\sum_{k=1}^{m}\tilde K_{1}^{(k)}(\nu ,\cdot )-\frac{1}{m}
\sum_{k=1}^{m}\tilde K_{2}^{(k)}(\nu ,\cdot )\right\Vert _{V}.  \notag
\end{eqnarray}%
By the definition of $k_{\varepsilon }$, each of the first two terms is
bounded by $\varepsilon +k_{\varepsilon }C_{\varepsilon}/m$, where $C_{\varepsilon}= \max\{\int V(w)\tilde K_i^{(k)}(\nu,dw): i = 1,2, ~~ 1\leq k\leq k_{\varepsilon} \}$.
For the last term, for each $k\leq m$ and $|f|\leq V$,
\begin{eqnarray*}
&&\left\vert \int f(w)[\tilde K_{1}^{(k+1)}(\nu ,dw)-\tilde K_{2}^{(k+1)}(\nu
,dw)]\right\vert  \\
&\leq &\left\vert \int \left( \tilde K_{1}^{(k)}(\nu ,dw)-\tilde K_{2}^{(k)}(\nu
,dw)\right) \int f(w')\tilde K_{2}(w,dw')\right\vert +\int \tilde K_{1}^{(k)}(\nu
,dw)\Vert \tilde K_{1}(w,\cdot)-K_{2}(w,\cdot)\Vert_V  \\
&\leq &\left\Vert \tilde K_{1}^{(k)}(\nu ,\cdot )-\tilde K_{2}^{(k)}(\nu ,\cdot
)\right\Vert _{V}+d(A_n) 
\\
&=&\left\Vert \tilde K_{1}^{(k)}(\nu ,\cdot )-\tilde K_{2}^{(k)}(\nu ,\cdot
)\right\Vert _{V} + o(1),
\end{eqnarray*}%
as $n\rightarrow \infty$.
The second inequality in the above display holds, if $V\equiv 1$; if $V$ is a drift function of $\tilde K_{2}$, we replace $f$ by $V$ and the inequality hold by noticing that $\int V(w')\tilde K_{2}(w,dw')\leq \lambda_{2} V(w)$.
Then, by induction, for all $k\leq m$,%
\begin{equation*}
\left\Vert \tilde K_{1}^{(k)}(\nu ,\cdot )-\tilde K_{2}^{(k)}(\nu ,\cdot )\right\Vert
_{V}\leq o(1)
\end{equation*}%
Therefore, the last term is
\begin{equation*}
\left\Vert \frac{1}{m}\sum_{k=1}^{m}\tilde K_{1}^{(k)}(\nu ,\cdot )-\frac{1}{m}%
\sum_{k=1}^{m}\tilde K_{2}^{(k)}(\nu ,\cdot )\right\Vert _{V}= o(1), 
\end{equation*}%
as $n\rightarrow \infty$.
Thus, for each $\varepsilon > 0$, we first choose $\kappa_{\varepsilon}$ and $C_{\varepsilon}$, then choose $m$ large such that $2C_{\varepsilon}k_{\varepsilon }/m < \varepsilon $, lastly choose $n$ large such that $\left\Vert \tilde K_{1}^{(k)}(\nu ,\cdot )-\tilde K_{2}^{(k)}(\nu ,\cdot )\right\Vert
_{V} < \varepsilon$.
Therefore,%
\begin{equation*}
\left\Vert \tilde \nu^{\xx^{obs}} _{1}-\tilde \nu^{\xx^{obs}} _{2}\right\Vert _{V}\leq
5\varepsilon .
\end{equation*}%
\end{proof}

\bigskip

\begin{proof}[Proof of Proposition \ref{PropDrift}]
The proof uses a similar idea as that of Proposition \ref{PropPrior}:  $\tilde K_1$ is equivalent to updating the missing values from the posterior predictive distribution of $f$ condition on that $\xx\in A_n$. Similarly, $\tilde K_2$ corresponds to the posterior predictive distributions of $g_j$'s.
By Proposition \ref{PropPrior},
for all $\xx\in A_{n}$, $|K_{1}(\xx,B) - K_{2}(\xx,B)|= O(n^{-1/4})$, which implies that
\begin{equation*}
\frac{K_{1}(w,A_{n})}{K_{2}(w,A_{n})}=1+O(n^{-1/4}).
\end{equation*}%
The posterior distribution is a joint distribution of the parameter and the missing values. Therefore, $\theta_{j}$ is part of the vector $w$. Let%
\begin{equation*}
R= \frac{\tilde K_2(w,dw')}{\tilde K_1(w,dw')}=\prod_{j=1}^{p}r_{j}(\theta _{j})\frac{K_{1}(w,A_{n})}{K_{2}(w,A_{n})},
\end{equation*}%
where $r_{j}(\theta_{j})$ is the normalized prior ratio corresponding to the imputation model of the $j$-th variable, whose definition is given in Proposition \ref{PropPrior}.


For the verification of the drift function, 
\begin{eqnarray}
\int V(w')\tilde{K}_{2}(w,dw^{\prime }) &=&\int R\times V(w^{\prime })\tilde{K}%
_{1}(w,dw^{\prime }) \notag  \\
&\leq &(1+O(n^{-1/4}))^{p}\int V(w^{\prime })\prod_{j=1}^{p}\left( 1+2\frac{%
\partial L(\mu _{\theta _{j}})Z_{j}}{L(\mu _{\theta _{j}})\sqrt{n}}\right) 
\tilde{K}_{1}(w,dw^{\prime }). \label{seq}
\end{eqnarray}%
Let $w_{j}$ be the state of the chain when the $j$-th variable is just
updated (then, $w' = w_{p}$). Then, according to the condition in \eqref{cond2}, we have  that for each $j+k\leq p$
\[
\tilde E_{1}\left[~ V(w_{j+k})\Big( 1+2\frac{\partial L(\mu _{\theta
_{j+k}})Z_{j+k}}{L(\mu _{\theta _{j+k}})\sqrt{n}}\Big)~ \Big |~ w_{j}~\right]
=\tilde E_{1}\left( V(w_{j+k})| w_{j}\right) +o(1)V(w_{j})
\]%
Since the $o(1)$ is uniform in $w_{j}\in A_{n}$, we can apply induction on
the product in (\ref{seq}) by conditioning on $\mathcal F_{j}=\sigma (w_{1},...,w_{j})$ sequentially for $j=1,...,n$. Therefore, we have 
\begin{eqnarray*}
\int V(w') \tilde{K}_{2}(w,dw') &=& (1+o(1))\int V(w')\tilde K_{1}(w,dw')+o(1)V(w) \\
&\leq &(\lambda _{1}+o(1))V(w).
\end{eqnarray*}%
Then, we can find another $\lambda _{2}\in (0,1)$ such that the above
display is bounded by $\lambda _{2}V(w)$. Thus, $V(w)$ is also a drift
function for $\tilde{K}_{2}$.
\end{proof}

\section{Proof of Theorem \ref{ThmIncomp}}

Throughout this proof we use the following notation for asymptotic
behavior. We say that $0\leq g(n)=O(h(n))$ if $g(n)\leq ch(n)$ for some
constant $c\in (0,\infty )$ and all $n\geq 1$.  We also write $g(n)=o(h(n))$ as $n\nearrow \infty $ if $g(n)/h(n)\rightarrow 0$ as $n\rightarrow \infty $. Finally, we write $X_{n} = O_{p}(g(n))$ if $|X_{n }/g(n)|$ is stochastically dominated by some distribution with finite exponential moment.

Let $\xx(k)$ be the iterative chain starting from its stationary distribution
$\nu _{2}^{\xx^{obs}}$. Furthermore, let $\nu _{2,j}^{\xx^{obs}}$ be the
distribution of $\xx(k)$ when the $j$-th variable is just updated.
 Due to incompatibility, $\nu _{2,j}^{\xx^{obs}}$'s are not necessarily identical.
Thanks to stationarity, $\nu _{2,j}^{\xx^{obs}}$ does not depend on $k$ and  $\nu _{2}^{\xx^{obs}}=\nu
_{2,0}^{\xx^{obs}}=\nu _{2,p}^{\xx^{obs}}$. Let
\begin{equation*}
(\tilde{\theta}_{j},\tilde{\varphi}_{j})=\arg \sup_{\theta_{j},\varphi_{j}} \int \log
h_{j}(\xx_{j}|\xx_{-j},\theta_j ,\varphi_j )\nu _{2,j-1}^{\xx^{obs}}(d\xx^{mis}).
\end{equation*}

The proof consists of two steps. Step 1, we show that for all $j$,
$\tilde{\varphi}_{j}\rightarrow 0,\tilde{\theta}_{j}-\hat{\theta}%
_{j}\rightarrow 0$, as $n\rightarrow \infty$, where $\hat \theta_{j}$ is the observed-data maximum likelihood estimate based on the joint model $f$ (defined as in \eqref{EE}).
That is, each variable is updated approximately from the conditional distribution $f(x_j|x_{-j},\hat \theta)$.
Step 2, we establish the statement of the theorem.

%

\paragraph{Step 1.}

We prove this step by contradiction. Suppose that there exist $\varepsilon
_{0}$ and $j_{0}$ such that $|\tilde{\varphi}_{j_{0}}|>\varepsilon _{0}$ or $%
|\tilde{\theta}_{j_{0}}-\hat{\theta}_{j_{0}}|>\varepsilon _{0}$. Let $\xx%
^{\ast }(k)$ be the Gibbs chain whose stationary distribution ($\nu _{1}^{\xx%
^{obs}}$) is the posterior predictive distribution associated with the joint model $f$ (c.f. Definition \ref{DefValid}). In addition, $\xx^{*}$ starts from its stationary distribution $\nu _{1}^{\xx%
^{obs}}$. We
now consider the KL divergence%
\[
D(\nu _{1}^{\xx^{obs}}||\nu _{2,j}^{\xx^{obs}}).
\]%
Let $\xx(k,j)$ and $\xx^{\ast }(k,j)$ be the state at iteration $k+1$ and
the $j$-th variable is just updated. 
Since both chains are stationary, the distributions of $\xx(k,j)$ and $\xx^{*}(k,j)$ are free of $k$.
To simplify notation, we let 
\begin{eqnarray*}
u &=&\xx_{j}^{mis}(0,j-1),\quad v=\xx_{-j}(0,j-1)=\xx_{-j}(0,j),\quad w=\xx%
_{j}^{mis}(0,j), \\
u^{\ast } &=&\xx_{j}^{\ast mis}(0,j-1),\quad v^{\ast }=\xx_{-j}^{\ast
}(0,j-1)=\xx_{-j}^{\ast }(0,j),\quad w^{\ast }=\xx_{j_{0}}^{\ast mis}(0,j).
\end{eqnarray*}%
That is, $u$ is the missing value of variable $j$ from the previous step and 
$w$ is the updated missing value of variable $j$. $v$ stands for the
variables that do not change in this update. Let $p_{j}(\cdot )$ be a
generic notation for density functions of $(u,v,w)$ and $p_{j}^{\ast }(\cdot )$
for $(u^{\ast },v^{\ast },w^{\ast })$. By the chain rule, we have that
\begin{eqnarray}
&&\int \log \frac{p_{j}^{\ast }(u,v,w)}{p_{j}(u,v,w)}p_{j}^{\ast
}(u,v,w)dudvdw \notag\\
&=&\int \log \frac{p_{j}^{\ast }(u,v)}{p_{j}(u,v)}%
p_{j}^{\ast }(u,v)dudv+\int \log \frac{p_{j}^{\ast }(w|u,v)}{p_{j}(w|u,v)}%
p_{j}^{\ast }(u,v,w)dudvdw \notag\\
&=&\int \log \frac{p_{j}^{\ast }(v,w)}{p_{j}(v,w)}p_{j}^{\ast
}(v,w)dvdw+\int \log \frac{p_{j}^{\ast }(u|v,w)}{p_{j}(u|v,w)}p_{j}^{\ast
}(u,v,w)dudvdw. \label{A}
\end{eqnarray}%
By construction,%
\begin{equation}\label{B}
\int \log \frac{p_{j}^{\ast }(u,v)}{p_{j}(u,v)}p_{j}^{\ast }(u,v)dudv=D(\nu
_{1}^{\xx^{obs}}||\nu _{2,j-1}^{\xx^{obs}})
\end{equation}
and%
\begin{equation}\label{C}
\int \log \frac{p_{j}^{\ast }(v,w)}{p_{j}(v,w)}p_{j}^{\ast }(v,w)dvdw=D(\nu
_{1}^{\xx^{obs}}||\nu _{2,j}^{\xx^{obs}}).
\end{equation}
Furthermore, $p_{j}^{\ast }(w|u,v)$ is the posterior predictive distribution
according to $f$ and $p_{j}(w|u,v)$ is the posterior predictive according to 
$h_{j}$. Note that $f$ is a sub-family of $h_{j}$ for the prediction of
variable $j$. 
Under $p_{j}^{\ast}(u,v,w)$ that is the posterior distribution of $f$, we have
$$\log \frac{p_{j}^{\ast }(w|u,v)}{f(\xx_{j}^{mis}|\xx_{-j},\hat \theta)}=O_{p}(1),\qquad \log \frac{p_{j}(w|u,v)}{f(\xx_{j}^{mis}|\xx_{-j},\hat \theta)}=O_{p}(1).$$
To understand the above estimate, for the posterior predictive distribution of the Gibbs chain, one first draw $\theta$ from the posterior distribution which is $\hat\theta+O_{p}(n^{-1/2})$ and then draw each $x_{j}$ from $f(x_{j}|x_{-j},\theta)$.
Thus, we may consider that
$p_{j}^{\ast }(w|u,v)\approx f(\xx_{j}^{mis}|\xx_{-j},\hat \theta + \xi/\sqrt n)$. Together with the fact that
$\partial \log f(\xx_{j}^{mis}|\xx_{-j},\hat \theta) = O_{p}(n^{-1/2})$  we obtain the above approximation. The same argument applies to the second estimate for $p_{j}(w|u,v)$ too.
With these two estimates, we have
$$\log \frac{p_{j}^{\ast }(w|u,v)}{p_{j}(w|u,v)}=O_{p}(1).$$
According to condition B3 that all the score functions has exponential moments, then we have
\begin{equation}\label{D}
\int \log \frac{p_{j}^{\ast }(w|u,v)}{p_{j}(w|u,v)}p_{j}^{\ast
}(u,v,w)dudvdw=O(1).
\end{equation}
We insert \eqref{B}, \eqref{C}, and \eqref{D} back to \eqref{A}.
For all $1\leq j\leq p$, we have that
\[
D(\nu _{1}^{\xx^{obs}}||\nu _{2,j-1}^{\xx^{obs}})=D(\nu _{1}^{\xx%
^{obs}}||\nu _{2,j}^{\xx^{obs}})+O(1)+\int \log \frac{p_{j}^{\ast }(u|v,w)}{%
p_{j}(u|v,w)}p_{j}^{\ast }(u,v,w)dudvdw.
\]%
We denote the last piece by 
\[
A_{j}=\int \log \frac{p_{j}^{\ast }(u|v,w)}{p_{j}(u|v,w)}p_{j}^{\ast
}(u,v,w)dudvdw.
\]%
 Note that $\nu _{2}^{\xx^{obs}}=\nu
_{2,0}^{\xx^{obs}}=\nu _{2,p}^{\xx^{obs}}$. Then, we have that 
\begin{eqnarray*}
D(\nu _{1}^{\xx^{obs}}||\nu _{2}^{\xx^{obs}}) &=&D(\nu _{1}^{\xx^{obs}}||\nu
_{2,0}^{\xx^{obs}}) \\
&=&D(\nu _{1}^{\xx^{obs}}||\nu _{2,1}^{\xx^{obs}})+A_{1}+O(1) \\
&=&... \\
&=&D(\nu _{1}^{\xx^{obs}}||\nu _{2,p}^{\xx^{obs}})+\sum_{j=1}^{p}A_{p}+O(1)
\\
&=&D(\nu _{1}^{\xx^{obs}}||\nu _{2}^{\xx^{obs}})+\sum_{j=1}^{p}A_{p}+O(1).
\end{eqnarray*}%
Thus, $\sum_{j=1}^{p}A_{p}=O(1)$. Note that each $A_{j}$ is non-negative.
Thus, $A_{j}$ must be bounded for all $j$, that is%
\begin{equation}
A_{j}=O(1).  \label{contra}
\end{equation}%
In what follows, we establish contradiction by showing that $A_{j_{0}}\rightarrow \infty $ if $|\tilde{\varphi}_{j_{0}}|+|\tilde{\theta}_{j_{0}}-\hat{\theta}_{j_{0}}|>\varepsilon _{0}$.

Now, we change all the $j$'s to $j_{0}$, that is, 
\begin{eqnarray*}
u &=&\xx_{j_{0}}^{mis}(0,j_{0}-1),\quad v=\xx_{-j_{0}}(0,j_{0}-1)=\xx%
_{-j_{0}}(0,j_{0}),\quad w=\xx_{j_{0}}^{mis}(0,j_{0}), \\
u^{\ast } &=&\xx_{j_{0}}^{\ast mis}(0,j_{0}-1),\quad v^{\ast }=\xx%
_{-j_{0}}^{\ast }(0,j_{0}-1)=\xx_{-j_{0}\ast }(0,j_{0}),\quad w^{\ast }=\xx%
_{j_{0}}^{\ast mis}(0,j_{0}).
\end{eqnarray*}%
Note that $u$ is the missing values of $\xx_{j_{0}}$ from the previous step
and $w$ is the missing value for the next step. In addition, the update of $%
\xx_{j_{0}}^{mis}$ does not depend on the previously imputed values.
Therefore, $u$ and $w$ are independent conditional on $v$. Thus, $A_{j_{0}}$
is reduced to%
\begin{eqnarray*}
\int \log \frac{p_{j_{0}}^{\ast }(u|v,w)}{p_{j_{0}}(u|v,w)}p_{j_{0}}^{\ast
}(u,v,w)dudvdw &=&\int \log \frac{p_{j_{0}}^{\ast }(u|v)}{p_{j_{0}}(u|v)}%
p_{j_{0}}^{\ast }(u,v)dudv \\
&=&\int \log \frac{d\nu _{1}^{\xx^{obs}}(\xx_{j_{0}}^{mis}|\xx_{-j_{0}})}{%
d\nu _{2,j_{0}-1}^{\xx^{obs}}(\xx_{j_{0}}^{mis}|\xx_{-j_{0}})}\nu _{1}^{\xx%
^{obs}}(d\xx^{mis}).
\end{eqnarray*}%
We further let $\iota $ be the set of observations where $x_{j_{0}}$ is
missing and $x_{-j_{0}}$ are observed. Use $\xx_{\iota ,j_{0}}^{mis}$ to
denote the missing $x_{j_{0}}$'s of the subset $\iota $. Then%
\[
\int \log \frac{p_{j_{0}}^{\ast }(u|v,w)}{p_{j_{0}}(u|v,w)}p_{j_{0}}^{\ast
}(u,v,w)dudvdw\geq \int \log \frac{d\nu _{1}^{\xx^{obs}}(\xx_{\iota
,j_{0}}^{mis}|\xx_{-j_{0}})}{d\nu _{2,j_{0}-1}^{\xx^{obs}}(\xx_{\iota
,j_{0}}^{mis}|\xx_{-j_{0}})}\nu _{1}^{\xx^{obs}}(d\xx^{mis}),
\]%
that is the joint K-L divergence is bounded from below by the marginal of
K-L divergence on the subset $\iota $. Note that $\xx_{\iota ,j_{0}}^{\ast
mis}(0,j_{0}-1)$ is the starting value of $\xx^{\ast }$ and was sampled from
the conditional stationary distribution $\nu _{1}^{\xx^{obs}}(\xx_{\iota
,j_{0}}^{mis}|\xx_{-j_{0}})$. Equivalently, $\xx_{\iota
,j_{0}}^{\ast mis}(0,j_{0}-1)$ is sampled from $f(\xx_{\iota ,j_{0}}|\xx%
_{\iota ,-j_{0}},\theta )$ where  $\theta = \hat \theta + O_{p}(n^{-1/2}) $ is a posterior sample  and $\xx_{\iota ,-j_{0}}$ is fully
observed (by the construction of set $\iota $).

On the other hand, $\xx_{\iota ,j_{0}}^{mis}(0,j_{0}-1)$ follows the
stationary distribution of the iterative chain and is sampled from the
previous step (step $k-1$) according to the conditional model $%
h_{j_{0}}(x_{j_{0}}|x_{i,-j_{0}},\theta _{j_{0}},\varphi _{j_{0}})$ where $%
(\theta _{j_{0}},\varphi _{j_{0}})= (\tilde{%
\theta}_{j_{0}},\tilde{\varphi}_{j_{0}}) +O_{p}(n^{-1/2})$ is a draw from the posterior distribution. In
addition, by assumption, the parameters are different by at least $%
\varepsilon _{0}$, that is, $|\tilde{\varphi}_{j_{0}}|+|\tilde{\theta}%
_{j_{0}}-\hat{\theta}_{j_{0}}|>\varepsilon _{0}$. Thus, the conditional distributions $h_{j_{0}}(\cdot |x_{i,-j_{0}},\tilde \theta_{j_{0}},\tilde \varphi_{j_{0}})$ and $f(\cdot |x_{i,-j_{0}},\hat \theta)$ are different. For some $\lambda _{0}>0$ (depending on $%
\varepsilon _{0}$), the KL divergence between the two updating distributions
of $x_{i,j_{0}}$ is bounded below by some $\lambda _{0}>0$, that is, for $i\in
\iota $ 
\[
D(~f(\cdot |x_{i,-j_{0}},\hat \theta)  ~\Vert ~ h_{j_{0}}(\cdot |x_{i,-j_{0}}, \tilde \theta_{j_{0}},\tilde \varphi_{j_{0}}) ~)\geq
\lambda _{0}.
\]%
This provides a lower bound of the KL divergence of one observation.
The posterior predictive distributions for the observations in $\iota $ are
conditionally independent given $(\theta _{j_{0}},\varphi _{j_{0}})$. Thus,
the KL divergence of the joint distributions is approximately the sum of the individual KL divergence of all the observations in $\iota $. Then,  we obtain that for some $\lambda_{1}>0$ 
\begin{eqnarray}
A_{j_{0}}=\int \log \frac{p_{j_{0}}^{\ast }(u|v,w)}{p_{j_{0}}(u|v,w)}p(u,v,w)dudvdw
&\geq &\int \log \frac{d\nu _{1}^{\xx^{obs}}(\xx_{\iota ,j_{0}}^{mis}|\xx%
_{-j_{0}})}{d\nu _{2,j_{0}-1}^{\xx^{obs}}(\xx_{\iota ,j_{0}}^{mis}|\xx%
_{-j_{0}})}\nu _{1}^{\xx^{obs}}(d\xx^{mis})  \label{far} \\
&\geq &\lambda _{1}\#(\iota ).  \nonumber
\end{eqnarray}%
%
%
Since the number of observations $\#(\iota )\rightarrow \infty $ (condition B5), we reached a
contradiction to (\ref{contra}). Thus, $|\tilde{\varphi}_{j}|+|\tilde{\theta}%
_{j}-\hat{\theta}_{j}|=o(1)$ as $n\rightarrow \infty $. Thereby, we conclude
Step 1.

\paragraph{Step 2.}

We first show the consistency of $\hat \theta^{(2)}$. It is sufficient to
show that $|\hat{\theta}^{(2)}-\hat{\theta}|\rightarrow 0$. $\hat{\theta}%
^{(2)}$ solves equation 
\[
\int \partial \log f(\xx^{mis},\xx^{obs}|\hat\theta^{(2)} )\nu _{2}^{\xx%
^{obs}}(d\xx^{mis})=0. 
\]%
By Taylor expansion, the MLE has the representation that 
\[
\hat \theta^{(2)} - \hat\theta = O(n^{-1})\int \partial \log f(\xx^{mis},\xx%
^{obs}~|~\hat\theta)\nu _{2}^{\xx^{obs}}(d\xx^{mis}). 
\]
Thus, it is sufficient to show that 
\[
\int \partial \log f(\xx^{mis},\xx^{obs}|\hat\theta)\nu _{2}^{\xx^{obs}}(d\xx%
^{mis})=o_{p}(n).
\]
Given that $\hat \theta - \theta^{0} = O_{p}(n^{-1/2})$. It is sufficient to
show that 
\begin{equation}  \label{cd}
\int \partial \log f(\xx^{mis},\xx^{obs}|\theta^{0})\nu _{2}^{\xx^{obs}}(d\xx%
^{mis})=o(n).
\end{equation}
Notice that the observed data MLE $\hat \theta$ satisfies 
\[
\int \partial \log f(\xx^{mis},\xx^{obs}|\hat \theta)f(d\xx^{mis}|\xx^{obs},%
\hat{\theta})=0 
\]%
and further 
\[
\int \partial \log f(\xx^{mis},\xx^{obs}|\theta^{0})f(d\xx^{mis}|\xx^{obs},%
\hat{\theta})=O_{p}(1). 
\]%
Then, we only need to show that 
\begin{equation}  \label{van}
\int \partial \log f(\xx^{mis},\xx^{obs}|\theta^{0})\Big [\nu _{2}^{\xx%
^{obs}}(d\xx^{mis})-f(d\xx^{mis}|\xx^{obs},\hat{\theta})\Big ] = o_{p}(n).
\end{equation}
Consider a single observation $x^{mis}(k)$. Without loss of generality,
suppose that $x^{mis}(k)=(x_{1}(k),...,x_{j}(k))$ and $%
x^{obs}=(x_{j+1},...,x_{p})$. 
The result of Step 1 suggests that each coordinate of $x^{mis}$ is updated
from 
\[
h_{j}(x_{j}|x_{-j},\hat \theta_{j}+o_{p}(1), \varphi_{j}= o_{p}(1)). 
\]
Thus, $x^{mis}(k)$ follows precisely the transition kernel of $x^*(k)$
described in condition B5. Therefore, we apply Lemma \ref{ThmComp} and have
that 
\begin{eqnarray*}
&&\int \partial \log f(x^{mis},x^{obs}|\theta^{0})\nu _{2}^{\xx%
^{obs}}(dx^{mis}) \\
&=&\int \partial \log f(x^{mis},x^{obs}|\theta^{0})f(x^{mis}|x^{obs},\hat{%
\theta})dx^{mis}+o(1).
\end{eqnarray*}%
Then, (\ref{van}) is satisfied immediately by adding up the above integral
for all observations. Therefore, \eqref{cd} is satisfied and further $\hat
\theta^{(2)} - \hat \theta \rightarrow 0$. The proof for $\hat \theta_j^{(2)}
$ and $\hat \varphi_j^{(2)}$ are completely analogous and therefore is
omitted. Thereby, we conclude the proof.

\section{Markov chain stability and rates of convergence}\label{SecMC}

In this section, we discuss the pending topic of the Markov chain's convergence. A bound on the convergence rate $q_k$ is required for both Lemma \ref{ThmComp} and \ref{ThmIncomp}. In this section, we review strategies in existing literature to check the convergence. We first provide a brief summary of methods to control the rate of convergence via renewal theory.

\paragraph{Markov chain stability by renewal theory.}
We first list a few conditions (cf.\ \cite{Bax05}), which we will refer to later.

\begin{enumerate}
  \item[C1]Minorization condition: A homogeneous Markov
process $W(n)$ with state space in $\mathcal X$ and
transition kernel $K(w,dw')= P(W(n+1) \in dw' | W(n) =w)$ is
said to satisfy a \emph{minorization condition} if for a subset
$C\subset \mathcal X$, there exists a probability measure $\nu$
on $\mathcal X$, $l\in \mathbb Z^+$, and $q\in (0,1]$ such
that
$$K^{(l)}(w,A) \geq q \nu(A)$$
for all $w\in C$ and measurable $A\subset\mathcal X$. $C$ is called a \emph{small set}.

  \item[C2] Strong aperiodicity condition: There exists $\delta >0$
  such that $q \nu(C) > \delta$.
  \item[C3] Geometric drift condition: there exists a non-negative and finite drift function,
$V$ and scalar $\lambda \in (0,1)$ such that for all $w\bar\in C$,
  $$\lambda V(w)\geq\int V(w') K(w,dw'),$$ and  for all
$w\in C$, $\int V(w') K(w,dw')\leq b$.
\end{enumerate}
Chains satisfying A1--3 are ergodic and admit a unique
stationary distribution
$$\pi(\cdot) =\lim_{n\rightarrow \infty} \frac 1 n \sum_{l=1}^n
K^{(l)}(w,\cdot)$$
for all $w$. Moreover, there exists $\rho<1$ depending only (and explicitly) on $q$, $\delta$, $\lambda$, and $b$ such that whenever $\rho <
\gamma <1$, there exists $M<\infty$ depending only (and
explicitly) on $q$, $\delta$, $\lambda$, and $b$ such that
\begin{equation}\label{cgtbd}\sup_{|g|\leq V} |\int g(w') K^{(k)}(w,dw')- \int g(w') \pi(dw')|\leq MV(w)\gamma^k,\end{equation}
for all $w$ and $k\geq 0$, where the supremum is
taken over all measurable $g$ satisfying $g(w) \leq V(w)$. See \cite{Ros95} and more recently \cite{Bax05} for a proof via the coupling of two Markov processes.

\paragraph{A practical alternative.}

In practice, one can check for convergence empirically. There are many diagnostic tools for the convergence of MCMC; see \cite{rhat} and the associated discussion.
Such empirical studies can show stability within the range of observed simulations.
This can be important in that we would like our imputations to be coherent even if
we cannot assure they are correct.
In addition, most theoretical bounds are conservative in the sense that the chain usually converges much faster than what it is implied by the bounds.  On the other hand, purely empirically checking supplies no theoretical guarantee that the chain converges to any distribution. Therefore, a theoretical development of the convergence is recommended when it is feasible given available resources (for instance, time constraint).




\bibliographystyle{plain}
\bibliography{bibstat}

\end{document}